\documentclass[10pt,a4paper]{article}
\usepackage{CJK}
\usepackage{mathrsfs}
\usepackage{amsfonts}
\begin{document}
\begin{CJK*}{GBK}{song}

\begin{center}
 {\LARGE\textbf{On Evaluation of Zeta and Related Functions by Abstract Operators}}
 \end{center}
\vspace{0.3cm}

\begin{center}
\rm Guang-Qing Bi
\end{center}

\begin{abstract}
Building on the mapping relations between analytic functions and periodic functions using the abstract operators $\cos(h\partial_x)$ and $\sin(h\partial_x)$, and by defining the Zeta and related functions including the Hurwitz Zeta function and the Dirichlet L-function in the form of abstract operators, we have obtained many new series expansions associated with these functions on the whole complex plane, and investigate the number theoretical properties of them, including some new rapidly converging series for $\eta(2n+1)$ and $\zeta(2n+1)$. For $n\in\mathbb{N}$, each of these series representing $\zeta(2n+1)$ converges remarkably rapidly with its general term having the order estimate:
\[O(m^{-2k}\cdot k^{-2n+1})\qquad(k\rightarrow\infty;\quad m=3,4,6).\]
\end{abstract}

\footnotetext{\hspace*{-.45cm} \noindent{\bf 2010 \emph{Mathematics Subject Classification}:}\ \ Primary 11M06, 35S05; Secondary 11B68, 11M35.}

\section{Preliminaries}
\noindent
\renewcommand\theequation{1.\arabic{equation}}

The Riemann Zeta function $\zeta(s)$, the Hurwitz Zeta function $\zeta(s,a)$ and the Dirichlet L-function $L(s,\chi)$ are defined usually by
\begin{equation}\label{00c}
\zeta(s):=\sum^\infty_{n=1}\frac{1}{n^s}\qquad(\Re(s)>1),
\end{equation}
\begin{equation}\label{00c'}
 \zeta(s,a):=\sum^\infty_{n=0}\frac{1}{(n+a)^s}\quad(\Re(s)>1,a\not\in\mathbb{Z}_0^-:=\{0,-1,-2,\ldots\})
\end{equation}
and for any Dirichlet character $\chi$ of modulus $q$ (There exists a positive integer $q$ such that $\chi(n)=\chi(n+q)$ for all $n$)
\begin{equation}\label{drlf}
 L(s,\chi):=\sum^\infty_{n=1}\frac{\chi(n)}{n^s}\quad(\Re(s)>1),
\end{equation}
and by their meromorphic continuations for $\Re(s)\leq1$. They are known to be meromorphic.

The alternating Zeta Function $\eta(s)$ is defined usually by
\begin{equation}\label{02}
\eta(s):=\sum^\infty_{n=1}\frac{(-1)^{n-1}}{n^s}=(1-2^{1-s})\zeta(s)\quad(\Re(s)>0),
\end{equation}
and by its analytic continuations for $\Re(s)\leq0$. It is a holomorphic function on the whole complex plane by analytic continuation.

The concept of abstract operators is based on the analytic continuity fundamental theorem contained in a 1997 paper entitled \emph{Applications of abstract operators to partial differential equations} (See \cite[p. 7]{bi97}). In other words, if two abstract operators $A$ and $B$ acts on $e^{\xi{x}},\,x\in\mathbb{R}^n,\,\xi\in\mathbb{R}^n$ such that $Ae^{\xi{x}}=Be^{\xi{x}}$, then $Ag(x)=Bg(x),\,\forall{g(x)}\in{C^\infty(\Omega)},\,\Omega\subset\mathbb{R}^n$. Therefore, the abstract operators $f(\partial_x),\,\partial_x:=(\partial_{x_1},\partial_{x_2},\cdots,\partial_{x_n})$ is defined by (See \cite[p. 8]{bi97})
\begin{equation}\label{abs}
 f(\partial_x)e^{\xi{x}}:=f(\xi)e^{\xi{x}},\quad\forall f(\xi)\in C^\infty(\mathbb{R}^n),
\end{equation}
where $f(\xi),\xi\in\mathbb{R}^n$ is called the symbols of abstract operators $f(\partial_x)$.

Further, we can derive the operation rules of abstract operators, such as
\begin{equation}\label{e0}
f(\partial_x)(e^{\xi_0x}g(x))=g(\partial_\xi)(e^{\xi{x}}f(\xi))|_{\xi=\xi_0},
\end{equation}
\begin{equation}\label{e1}
f(\partial_x)(e^{\xi{x}}g(x))=e^{\xi{x}}f(\xi+\partial_x)g(x),
\end{equation}
\begin{equation}\label{e1'}
 e^{h\partial_x}f(x)=f(x+h),\quad h\partial_x:=h_1\partial_{x_1}+h_2\partial_{x_2}+\cdots+h_n\partial_{x_n}
\end{equation}
and so on. The abstract operators in view of the analytic continuity fundamental theorem can be called the pseudo-differential operators defined on $C^\infty(\Omega)$. Conversely, the pseudo-differential operators in view of the Fourier transform can also be called the abstract operators defined on $\mathscr{S}(\mathbb{R}^n)$.

The abstract operators have been applied strongly to partial differential equations (See \cite{bi18}-\cite{bi11}). In this paper, we will see that the concept of abstract operators can be applied in defining alternating Zeta function, Riemann Zeta function, Hurwitz Zeta function and Dirichlet L-function on the whole complex plane without the experience of analytic continuation process. From this, we can easily derive their asymptotic expansions and other arithmetic properties, including some new rapidly converging series for $\eta(2n+1)$ and $\zeta(2n+1)$.

Let $h\partial_x=\langle{h,\partial_x}\rangle=h_1\partial_{x_1}+h_2\partial_{x_2}+\cdots+h_n\partial_{x_n}$. Then $\cos(h\partial_x)$ and $\sin(h\partial_x)$ are the abstract operators taking $\cos(hb)$ and $\sin(hb)$ as the symbols respectively, namely
\begin{equation}\label{1}
    \cos(h\partial_x)e^{bx}:=\cos(bh)e^{bx},\quad\sin(h\partial_x)e^{bx}:=\sin(bh)e^{bx}.
\end{equation}
Here $bx=b_1x_1+b_2x_2+\cdots+b_nx_n,\;hb=bh=b_1h_1+b_2h_2+\cdots+b_nh_n$. Further, their operation rules can be expressed as the following three groups of operator relationships by Guang-Qing Bi \cite[pp. 7-9]{bi97}:

\textbf{Theorem 1.1.} (See \cite[p. 9, Theorem 3, 4, and 6]{bi97}) Let $x\in\mathbb{R}^n,\;h\in\mathbb{R}^n$, $h\partial_x=\langle{h,\partial_x}\rangle=h_1\partial_{x_1}+\cdots+h_n\partial_{x_n}$. For $\cos(h\partial_x)\,\mbox{and}\,\sin(h\partial_x)$, we have
\begin{equation}\label{y0}
  \cos(h\partial_x)f(x)=\Re[f(x+ih)],\quad\sin(h\partial_x)f(x)=\Im[f(x+ih)],
\end{equation}
$\forall{f(z)}\in{C}^\infty(\Omega),\,z=x+iy\in\Omega\subseteq\mathbb{C}^n$;

\parbox{11cm}{\begin{eqnarray*}\label{y1}
                \sin(h\partial_x)(uv) &=& \cos(h\partial_x)v\cdot\sin(h\partial_x)u+\sin(h\partial_x)v\cdot\cos(h\partial_x)u,\\
                \cos(h\partial_x)(uv) &=& \cos(h\partial_x)v\cdot\cos(h\partial_x)u-\sin(h\partial_x)v\cdot\sin(h\partial_x)u;
              \end{eqnarray*}}\hfill\parbox{1cm}{\begin{eqnarray}\end{eqnarray}}

\parbox{11cm}{\begin{eqnarray*}\label{y2}
                \sin(h\partial_x)\frac{u}{v} &=& \frac{\cos(h\partial_x)v\cdot\sin(h\partial_x)u-\sin(h\partial_x)v\cdot\cos(h\partial_x)u}
                {(\cos(h\partial_x)v)^2+(\sin(h\partial_x)v)^2}, \\
                \cos(h\partial_x)\frac{u}{v} &=& \frac{\cos(h\partial_x)v\cdot\cos(h\partial_x)u+\sin(h\partial_x)v\cdot\sin(h\partial_x)u}
                {(\cos(h\partial_x)v)^2+(\sin(h\partial_x)v)^2}.
              \end{eqnarray*}}\hfill\parbox{1cm}{\begin{eqnarray}\end{eqnarray}}

\textbf{Theorem 1.2.} (See \cite[p. 9, Theorem 5]{bi97}) Let $h_0\in\mathbb{R},\;x(t)\in\mathbb{R}^n,\;t\in\mathbb{R}^1,\;X\in\mathbb{R}^n,\,Y\in\mathbb{R}^n$,
$Y\partial_X=\langle{Y,\partial_X}\rangle=Y_1\partial_{X_1}+\cdots+Y_n\partial_{X_n}$.
Then we have

\parbox{11cm}{\begin{eqnarray*}\label{y3}
 \sin(h_0\partial_t)f(x(t)) &=& \sin(Y\partial_X)f(X),\\
 \cos(h_0\partial_t)f(x(t)) &=& \cos(Y\partial_X)f(X),
\end{eqnarray*}}\hfill\parbox{1cm}{\begin{eqnarray}\end{eqnarray}}
where $X_j=\cos(h_0\partial_t)x_j(t),\;Y_j=\sin(h_0\partial_t)x_j(t),\;j=1,\cdots,n$.

In the special case when $n=1$, (\ref{y3}) can easily be restated as

\parbox{11cm}{\begin{eqnarray*}\label{yb1}
\sin\left(h_0\frac{d}{dt}\right)f(x(t)) &=& \sin\left(Y\frac{\partial}{\partial{X}}\right)f(X),\\
\cos\left(h_0\frac{d}{dt}\right)f(x(t)) &=& \cos\left(Y\frac{\partial}{\partial{X}}\right)f(X),
\end{eqnarray*}}\hfill\parbox{1cm}{\begin{eqnarray}\end{eqnarray}}
where
$Y=\sin(h_0\partial_t)x(t),\;X=\cos(h_0\partial_t)x(t),\;t\in\mathbb{R}^1,\;h_0\in\mathbb{R}$.

\textbf{Theorem 1.3.} (See \cite[p. 9, Theorem 7]{bi97}) Let $u=g(y)$ be a monotonic function on its domain. If $y=f(bx)$ is the inverse function of $bx=g(y)$ such that $g(f(bx))=bx$, where $bx=b_1x_1+b_2x_2+\cdots+b_nx_n,\;bh=b_1h_1+b_2h_2+\cdots+b_nh_n$,
then $\sin(h\partial_x)f(bx)$ (denoted by $Y$) and $\cos(h\partial_x)f(bx)$ (denoted by $X$) can be determined by the following set of equations:
\begin{equation}\label{y4}
    \left\{\begin{array}{l@{\qquad}l}\displaystyle
    \cos\left(Y\frac{\partial}{\partial{X}}\right)g(X)=bx,&x\in\mathbb{R}^n,\;b\in\mathbb{R}_n,\\\displaystyle
    \sin\left(Y\frac{\partial}{\partial{X}}\right)g(X)=bh, &h\in\mathbb{R}^n.
    \end{array}\right.
\end{equation}

For example, by making use of (\ref{1}) and (\ref{y4}), we have
\begin{eqnarray*}
    \left\{\begin{array}{l@{\qquad}l}\displaystyle
    e^X\cos{Y} = bx, & X=\cos(h\partial_x)\ln(bx),\\\displaystyle
    e^X\sin{Y} = bh, & Y=\sin(h\partial_x)\ln(bx).
    \end{array}\right.
\end{eqnarray*}
By solving this set of equations, we obtain

\parbox{10cm}{\begin{eqnarray*}\label{7}
                \cos(h\partial_x)\ln(bx) &=& \frac{1}{2}\ln\left((bx)^2+(bh)^2\right), \\
                \sin(h\partial_x)\ln(bx) &=&\textrm{arccot}\frac{bx}{bh}.
              \end{eqnarray*}}\hfill\parbox{1cm}{\begin{eqnarray}\end{eqnarray}}

\textbf{Theorem 1.4.}  Let $f(x)\in{L^2}[-c,c]$ be the sum function of the Fourier cosine series, and $g(x)\in{L^2}[-c,c]$ be that of the corresponding Fourier  sine series, namely
\[ f(x)=\sum^\infty_{n=0}a_n\cos\frac{n\pi{x}}{c}\quad\mbox{and}\quad g(x)=\sum^\infty_{n=0}a_n\sin\frac{n\pi{x}}{c},\]
where $x\in\Omega\subseteq\mathbb{R}^1$, $c>0$ is a given real number. If $S(t)$ is the sum function of the corresponding power series $\sum^\infty_{n=0}a_nt^n$, namely
\[S(t)=\sum^\infty_{n=0}a_nt^n,\quad(t\in\mathbb{R}^1,\;|t|<r,\;0<r<+\infty),\]
then for $a<x<b$ we have the following mapping relationships:

\parbox{11cm}{\begin{eqnarray*}\label{12}
                f(x) &=& \left.\cos\left(\frac{\pi{x}}{c}\frac{\partial}{\partial{z}}\right)S(e^z)\right|_{z=0}, \\
                g(x) &=& \left.\sin\left(\frac{\pi{x}}{c}\frac{\partial}{\partial{z}}\right)S(e^z)\right|_{z=0}.
              \end{eqnarray*}}\hfill\parbox{1cm}{\begin{eqnarray}\end{eqnarray}}
Here interval $(a,b)\subset\Omega$. The endpoints $a$ and $b$ of the interval $a<x<b$ are non-analytical points (singularities) of Fourier series, which can be uniquely determined by the detailed computation of the right-hand side of (\ref{12}).

\textbf{Proof.} Theorem 1.4 can be proved easily by substituting $S(e^z)=\sum^\infty_{n=0}a_ne^{nz}$ into (\ref{12}).

In other words, (\ref{12}) gives the following trigonometric summation relationships:
\begin{equation}\label{12ct1}
 \left.\cos\left(\frac{\pi{x}}{c}\frac{\partial}{\partial{z}}\right)S(e^z)\right|_{z=0}=\sum^\infty_{n=0}a_n\cos\frac{n\pi{x}}{c};
\end{equation}
\begin{equation}\label{12ct2}
 \left.\sin\left(\frac{\pi{x}}{c}\frac{\partial}{\partial{z}}\right)S(e^z)\right|_{z=0}=\sum^\infty_{n=0}a_n\sin\frac{n\pi{x}}{c}.
\end{equation}
Here $x\in\Omega\subseteq\mathbb{R}^1$. The $\Omega$ can be uniquely determined by the detailed computation of the left-hand side of (\ref{12ct1}) and (\ref{12ct2}) respectively.

For example, according to the proof of Lemma 2.3 in second sections of this paper, if $S(e^z)=\ln(1+e^z)$, then we have $\Omega:=\{x\in\mathbb{R}^1|\cos(\pi{x}/(2c))>0\}$ for (\ref{12ct1}) and $\Omega:=\{x\in\mathbb{R}^1|\cos(\pi{x}/(2c))\neq0\}$ for (\ref{12ct2}) respectively.

\textbf{Theorem 1.5.}  Let $S(t)$ be an arbitrary analytic function integrable in the interval $[0,1]$. Then we have
\begin{equation}\label{13}
  \left.\cos\left(\frac{\pi{x}}{c}\frac{\partial}{\partial{z}}\right)\int^{e^z}_0\!\!\!S(e^z)\,de^z\right|_{z=0}
 = \int^1_0\!\!S(\xi)\,d\xi-\frac{\pi}{c}\int^x_0\!\!\left.\sin\left(\frac{\pi{x}}{c}\frac{\partial}{\partial{z}}\right)
 [S(e^z)\,e^z]\right|_{z=0}dx.
\end{equation}
\begin{equation}\label{14}
 \left.\sin\left(\frac{\pi{x}}{c}\frac{\partial}{\partial{z}}\right)\int^{e^z}_0\!\!S(e^z)\,de^z\right|_{z=0}=
\frac{\pi}{c}\int^x_0\!\!\left.\cos\left(\frac{\pi{x}}{c}\frac{\partial}{\partial{z}}\right)[S(e^z)\,e^z]\right|_{z=0}dx.
\end{equation}

\textbf{Proof.} According to the analytic continuity fundamental theorem, we only need to prove this set of formulas when $S(x)=x^n$, $n\in\mathbb{N}_0:=\{0,1,2,\cdots\}$. This is obvious.

Let $S_0(t)$ be a function analytic in the neighborhood of $t=0$ and
\[S_0(t)=\sum^\infty_{n=1}a_nt^n,\qquad t\in\mathbb{R}^1,\;\;|t|<r,\quad0<r<+\infty,\]
where $a_n$ are rational numbers, then the sum function $S_m(t)$ is defined as
\begin{equation}\label{23}
S_m(t):=\underbrace{\int^t_0\frac{dt}{t}\cdots}_m\int^t_0S_0(t)\,\frac{dt}{t}=\sum^\infty_{n=1}a_n\frac{t^n}{n^m},\quad m\in\mathbb{N}_0:=\{0,1,2,\cdots\}.
\end{equation}
Apparently $S_m(1)$ is the sum function of the Dirichlet series taking $m$ as the variable.

According to (\ref{23}), $S_m(t)$ satisfies the following recurrence relation:
\begin{equation}\label{24}
\int^t_0S_{m-1}(t)\,\frac{dt}{t}=S_m(t),\qquad m\in\mathbb{N}:=\{1,2,\cdots\}.
\end{equation}

\textbf{Theorem 1.6.}  Let $m\in\mathbb{N},\,c>0$. The sum function $S_m(t)$ has the following recurrence property:
\begin{equation}\label{25}
\left.\cos\left(\frac{\pi{x}}{c}\frac{\partial}{\partial{z}}\right)S_m(e^z)\right|_{z=0}=S_m(1)
-\frac{\pi}{c}\int^x_0\left.\sin\left(\frac{\pi{x}}{c}\frac{\partial}{\partial{z}}\right)S_{m-1}(e^z)\right|_{z=0}dx,
\end{equation}
\begin{equation}\label{26}
\left.\sin\left(\frac{\pi{x}}{c}\frac{\partial}{\partial{z}}\right)S_m(e^z)\right|_{z=0}=
\frac{\pi}{c}\int^x_0\left.\cos\left(\frac{\pi{x}}{c}\frac{\partial}{\partial{z}}\right)S_{m-1}(e^z)\right|_{z=0}dx.
\end{equation}

\textbf{Proof.}  Taking $S(x)=S_{m-1}(x)/x$ in Theorem 1.5, then it is proved by using (\ref{24}).

In particular, if $S_0(t)=\ln(1+t)$, then $S_m(1)=\eta(m+1)$; if $S_0(t)=-\ln(1-t)$, then $S_m(1)=\zeta(m+1)$.

\section{Lemmas and definitions}
\noindent\setcounter{equation}{0}
\renewcommand\theequation{2.\arabic{equation}}

\textbf{Lemma 2.1.} For $m,r\in\mathbb{N}$, the sum function $S_m(t)$ has the following recurrence property:
\begin{eqnarray}\label{30"}
                & & \left.\cos\left(\frac{\pi{x}}{c}\frac{\partial}{\partial{z}}\right)S_{m-1}(e^z)\right|_{z=0}\nonumber\\
                &=& \sum^{r-1}_{k=0}(-1)^k\frac{1}{(2k)!}\left(\frac{\pi{x}}{c}\right)^{2k}S_{m-2k-1}(1)\nonumber\\
                & & +\,(-1)^r\left(\frac{\pi}{c}\right)^{2r-1}
\underbrace{\int^x_0dx\cdots}_{2r-1}\int^x_0\left.\sin\left(\frac{\pi{x}}{c}\frac{\partial}{\partial{z}}\right)S_{m-2r}(e^z)\right|_{z=0}dx.
\end{eqnarray}

\textbf{Proof.}  We can use the mathematical induction to prove it. According to (\ref{25}) of Theorem 1.6, it is obviously tenable when $r=1$ in (\ref{30"}). Now we inductively hypothesize that it is tenable when $r=K$. Using Theorem 1.6, then it is tenable when $r=K+1$, thus Lemma 2.1 is proved.

Similarly,

\textbf{Lemma 2.2.} For $m,r\in\mathbb{N}$, the sum function $S_m(t)$ has the following recurrence property:

\begin{eqnarray}\label{30}
                & & \left.\cos\left(\frac{\pi{x}}{c}\frac{\partial}{\partial{z}}\right)S_{m-1}(e^z)\right|_{z=0}\nonumber\\
                &=& \sum^{r-1}_{k=0}(-1)^k\frac{1}{(2k)!}\left(\frac{\pi{x}}{c}\right)^{2k}S_{m-2k-1}(1)\nonumber\\
                & & +\,(-1)^r\left(\frac{\pi}{c}\right)^{2r}
\underbrace{\int^x_0dx\cdots}_{2r}\int^x_0\left.\cos\left(\frac{\pi{x}}{c}\frac{\partial}{\partial{z}}\right)S_{m-2r-1}(e^z)\right|_{z=0}dx.
\end{eqnarray}

\textbf{Lemma 2.3.}  Let $x\in\mathbb{R}^1$ with $|x|<c$. Then

\begin{equation}\label{yc7.3+}
\left.\sin\left(\frac{\pi{x}}{c}\frac{\partial}{\partial{z}}\right)\ln(1+e^z)\right|_{z=0}=\frac{\pi{x}}{2c}.
 \end{equation}
\begin{equation}\label{yc7.3}
\left.\cos\left(\frac{\pi x}{c}\frac{\partial}{\partial z}\right)\ln(1+e^z)\right|_{z=0}=\ln\left(2\cos\frac{\pi x}{2c}\right).
\end{equation}

\textbf{Proof.} By making use of (\ref{1}), (\ref{yb1}) and (\ref{7}), we have
\begin{eqnarray*}
   & &\left.\sin\left(\frac{\pi{x}}{c}\frac{\partial}{\partial{z}}\right)\ln(1+e^z)\right|_{z=0}\,=\,
   \left.\sin\left(Y\frac{\partial}{\partial{X}}\right)\ln{X}\right|_{z=0}=\left.\textrm{arccot}\frac{X}{Y}\right|_{z=0}\\
   &=&\textrm{arccot}\frac{1+\cos(\pi{x}/c)}{\sin(\pi{x}/c)}=\textrm{arccot}\frac{\cos^2(\pi x/(2c))}{\sin(\pi x/(2c))\cos(\pi x/(2c))}.
\end{eqnarray*}
When $\cos(\pi{x}/(2c))\neq0$ or $|x|\neq{c}$, the above expression can be written in the form:
\[\left.\sin\left(\frac{\pi{x}}{c}\frac{\partial}{\partial{z}}\right)\ln(1+e^z)\right|_{z=0}=
\textrm{arccot}\cot\frac{\pi{x}}{2c}=\frac{\pi{x}}{2c}\quad(|x|<c).\]

Similarly, for $\cos(\pi{x}/(2c))>0$, we have
\begin{eqnarray*}
   & & \left.\cos\left(\frac{\pi{x}}{c}\frac{\partial}{\partial{z}}\right)\ln(1+e^z)\right|_{z=0}=
   \left.\cos\left(Y\frac{\partial}{\partial{X}}\right)\ln{X}\right|_{z=0}=\left.\frac{1}{2}\ln(X^2+Y^2)\right|_{z=0}\\
   &=& \frac{1}{2}\ln\left(\left(1+\cos\frac{\pi{x}}{c}\right)^2+\sin^2\frac{\pi{x}}{c}\right)=\ln\left(2\cos\frac{\pi{x}}{2c}\right)\quad(|x|<c).
\end{eqnarray*}
Thus Lemma 2.3 is proved.

\textbf{Lemma 2.4.} For $\Re(s)>1$
\begin{equation}\label{yc2}
\sum^\infty_{n=1}\frac{1}{n^s}\cos\frac{2n\pi}{3}=\frac{1}{2}(3^{1-s}-1)\zeta(s).
\end{equation}
\begin{equation}\label{yc4}
\sum^\infty_{n=1}\frac{1}{n^s}\cos\frac{n\pi}{2}=2^{-s}(2^{1-s}-1)\zeta(s).
\end{equation}
\begin{equation}\label{yc1}
\sum^\infty_{n=1}\frac{1}{n^s}\cos\frac{n\pi}{3}=\frac{1}{2}(6^{1-s}-3^{1-s}-2^{1-s}+1)\zeta(s).
\end{equation}

\textbf{Proof.} Let $\Phi(z,s,a)$ be the Lerch transcendent.  For $\Re(s)>1$, it is given by
\begin{equation}\label{ycz0}
 \Phi(z,s,a)=\sum^\infty_{n=0}\frac{z^n}{(n+a)^s}\quad(a\not\in\mathbb{Z}_0^-,\,|z|\leq1),
\end{equation}
so that $\Phi(1,s,a)=\zeta(s,a)$ and
\begin{equation}\label{ycz1}
 \zeta(s,a)+\Phi(-1,s,a)=2^{1-s}\zeta(s,a/2).
\end{equation}
Since $\zeta(s,1)=\zeta(s),\;\zeta(s,2)=\zeta(s)-1$ and
\begin{equation}\label{ycz2}
 \sum^{m-1}_{k=0}\zeta(s,a+k/m)=m^s\zeta(s,ma)\quad(m\in\mathbb{N}),
\end{equation}
we obtain
\begin{equation}\label{ycz3}
 \Phi(-1,s,2/3)+\Phi(-1,s,4/3)=(2^{1-s}-1)(3^s-1)\zeta(s)+3^s.
\end{equation}
Thus we have
\begin{eqnarray*}
   \sum^\infty_{n=1}\frac{1}{n^s}\cos\frac{n\pi}{3} &=& \cos\frac{\pi}{3}+\sum^\infty_{n=1}\frac{1}{(3n-1)^s}\cos\left(n\pi-\frac{\pi}{3}\right) \\
   & & +\sum^\infty_{n=1}\frac{1}{(3n)^s}\cos(n\pi)+\sum^\infty_{n=1}\frac{1}{(3n+1)^s}\cos\left(n\pi+\frac{\pi}{3}\right)\\
   &=& \frac{3^s-(\Phi(-1,s,2/3)+\Phi(-1,s,4/3))}{2\times3^s}-\frac{\eta(s)}{3^s}\\
   &=& -\frac{(2^{1-s}-1)(3^s-1)\zeta(s)}{2\times3^s}-\frac{(1-2^{1-s})\zeta(s)}{3^s}\\
   &=& \frac{1}{2}(6^{1-s}-3^{1-s}-2^{1-s}+1)\zeta(s)\qquad(\Re(s)>1).
\end{eqnarray*}
Similarly we have (\ref{yc2}) and (\ref{yc4}). Lemma 2.4 is proved.

\textbf{Lemma 2.5.} For $\Re(s)>1$
\begin{equation}\label{yc2.7}
\sum^\infty_{n=1}(-1)^{n-1}\frac{1}{n^s}\cos\frac{2n\pi}{3}=-\frac{6^{1-s}-3^{1-s}-2^{1-s}+1}{2(1-2^{1-s})}\eta(s).
\end{equation}
\begin{equation}\label{yc4.7}
\sum^\infty_{n=1}(-1)^{n-1}\frac{1}{n^s}\cos\frac{n\pi}{2}=2^{-s}(1-2^{1-s})\zeta(s)=2^{-s}\eta(s).
\end{equation}
\begin{equation}\label{yc1.7}
\sum^\infty_{n=1}(-1)^{n-1}\frac{1}{n^s}\cos\frac{n\pi}{3}=\frac{1-3^{1-s}}{2(1-2^{1-s})}\eta(s).
\end{equation}

\textbf{Proof.} It is easily seen that for $0\leq x\leq c$
\begin{equation}\label{yc7.1}
\sum^\infty_{n=1}(-1)^{n-1}\frac{1}{n^s}\cos\frac{n\pi(c-x)}{c}=-\sum^\infty_{n=1}\frac{1}{n^s}\cos\frac{n\pi x}{c},\quad(\Re(s)>1).
\end{equation}
By applying relationships (\ref{yc7.1}) in the case $x=c/3, c/2$ and $2c/3$ respectively, and by using the relation $\eta(s)=(1-2^{1-s})\zeta(s)$, we obtain Lemma 2.5 from Lemma 2.4.

\textbf{Lemma 2.6.} Let $x\in\mathbb{R}^1$ with $|x|<c$. Then

\begin{equation}\label{yc7.2}
\ln\cos\frac{\pi x}{2c}=\sum^\infty_{k=1}(-1)^k\frac{E_{2k-1}(1)}{2(2k)!}\left(\frac{\pi x}{c}\right)^{2k},
\end{equation}
where $E_n(x)$ are the Euler polynomials defined by the generating functions:
\begin{equation}\label{sch}
  \frac{2e^{xz}}{e^z+1}=\sum^\infty_{n=0}E_n(x)\frac{z^n}{n!}\quad(|z|<\pi).
\end{equation}

\textbf{Proof.} By (\ref{yc7.3}) of Lemma 2.3, we have
\begin{eqnarray*}
 \ln\cos\frac{\pi x}{2c}&=&\left.\cos\left(\frac{\pi{x}}{c}\frac{\partial}{\partial{z}}\right)\ln(1+e^z)\right|_{z=0}-\ln2\\
   &=& \sum^\infty_{k=0}\frac{(-1)^k}{(2k)!}\left(\frac{\pi x}{c}\right)^{2k}\!\left.\frac{\partial^{2k}}{\partial z^{2k}}\ln(1+e^z)\right|_{z=0}-\ln2\\
   &=& \sum^\infty_{k=1}\frac{(-1)^k}{(2k)!}\left(\frac{\pi x}{c}\right)^{2k}\!\left.\frac{\partial^{2k-1}}{\partial z^{2k-1}}\frac{e^z}{1+e^z}\right|_{z=0}\\
   &=& \sum^\infty_{k=1}\frac{(-1)^k}{(2k)!}\left(\frac{\pi x}{c}\right)^{2k}\frac{E_{2k-1}(1)}{2}\quad(|x|<c).
\end{eqnarray*}

\textbf{Lemma 2.7.} For $r\in\mathbb{N}$ and $|x|\leq c\;(x\in\mathbb{R}^1)$,
\begin{eqnarray}\label{d7.1}
 \sum^\infty_{n=1}\frac{(-1)^{n-1}}{n^{2r+1}}\cos\frac{n\pi{x}}{c}
 &=& \sum^r_{k=0}\frac{(-1)^k}{(2k)!}\left(\frac{\pi{x}}{c}\right)^{2k}\eta(2r+1-2k)\nonumber\\
 & & +\sum^\infty_{k=1}(-1)^{r+k}\frac{E_{2k-1}(1)}{2(2r+2k)!}\left(\frac{\pi{x}}{c}\right)^{2r+2k},
\end{eqnarray}
which is a slightly corrected version of a result proven in a significantly different way by Tsumura \cite[p. 388, Proposition 1 (3)]{Tsu}.

\textbf{Proof.} In Lemma 2.2, let $m=2r+1,\;S_0(t)=\ln(1+t)$. By using Theorem 1.4 and Lemma 2.3, since
$S_{2r-2k}(1)=\eta(2r+1-2k)$,  we have
\begin{eqnarray}\label{d7.1"}
 \sum^\infty_{n=1}\frac{(-1)^{n-1}}{n^{2r+1}}\cos\frac{n\pi{x}}{c}
 &=& \sum^{r-1}_{k=0}\frac{(-1)^k}{(2k)!}\left(\frac{\pi{x}}{c}\right)^{2k}\eta(2r+1-2k)\nonumber\\
 & & +\,(-1)^r\left(\frac{\pi}{c}\right)^{2r}\underbrace{\int^x_0dx\cdots}_{2r}\int^x_0\ln\left(2\cos\frac{\pi x}{2c}\right)dx,
\end{eqnarray}
where $|x|\leq c\;(x\in\mathbb{R}^1)$ for $r\in\mathbb{N}$. Substituting (\ref{yc7.2}) into (\ref{d7.1"}), since $\eta(1)=\ln2$, then it is proved.

\textbf{Lemma 2.8.} For $r\in\mathbb{N}$ and $|x|\leq c\;(x\in\mathbb{R}^1)$,
\begin{eqnarray}\label{44L}
 \sum^{\infty}_{n=1}(-1)^{n-1}\frac{1}{n^{2r}}\cos\frac{n\pi{x}}{c}
 =\sum^r_{k=0}(-1)^{k}\frac{1}{(2k)!}\left(\frac{\pi{x}}{c}\right)^{2k}\eta(2r-2k),
\end{eqnarray}
which is a slightly corrected version of a result proven in a significantly different way by Tsumura \cite[p. 387, Proposition 1 (1)]{Tsu}.

\textbf{Proof.} Letting $m=2r$ in (\ref{30"}), if $S_0(t)=\ln(1+t)$, then according to Lemma 2.1, Theorem 1.4 and Lemma 2.3, we can obtain (on the interval $[-c,c]$ for $r\in\mathbb{N}$):
\begin{eqnarray}\label{44}
 \sum^{\infty}_{n=1}(-1)^{n-1}\frac{1}{n^{2r}}\cos\frac{n\pi{x}}{c}
 &=&\sum^{r-1}_{k=0}(-1)^{k}\frac{1}{(2k)!}\left(\frac{\pi{x}}{c}\right)^{2k}\eta(2r-2k)\nonumber\\
 & &+(-1)^r\frac{1}{2(2r)!}\left(\frac{\pi{x}}{c}\right)^{2r}.
\end{eqnarray}
Since $\eta(0)=1/2$, then it is proved.

Let $\partial_z^{-s},\,z\in\mathbb{R}^1$ be an abstract operators taking $\xi^{-s}$ as the symbols, namely
\begin{equation}\label{cxsz}
 \left(\frac{\partial}{\partial z}\right)^{-s}e^{\xi z}:=\xi^{-s}e^{\xi z}\quad(\xi\in\mathbb{R}^1,\,s\in\mathbb{C}).
\end{equation}

\textbf{Definition 2.1.} The alternating Zeta Function $\eta(s)$ can be defined in the complete form:
\begin{equation}\label{yrdy}
 \eta(s):=\left(\frac{\partial}{\partial z}\right)^{-s}\!\!\left.\frac{e^z}{1+e^z}\right|_{z=0},\quad s\in\mathbb{C},
\end{equation}
which is a holomorphic function on the whole complex plane.

Similarly,
\begin{equation}\label{yrkz}
 \eta(s)=\sum^m_{n=1}\frac{(-1)^{n-1}}{n^s}+(-1)^m\left(\frac{\partial}{\partial z}\right)^{-s}\!\!\left.\frac{e^{(m+1)z}}{1+e^z}\right|_{z=0}\quad (m\in\mathbb{N}_0).
\end{equation}

\textbf{Definition 2.2.} The Hurwitz Zeta function $\zeta(s,a)$ can be defined in the complete form:
\begin{equation}\label{hurw}
 \zeta(s,a):=\frac{1}{s-1}\left(\frac{\partial}{\partial z}\right)^{1-s}\!\!\left.\frac{ze^{az}}{e^z-1}\right|_{z=0}\quad(a\not\in\mathbb{Z}_0^-,\,s\in\mathbb{C}\;\mbox{and}\;s\neq1).
\end{equation}

This definition is valid for all complex $s$. In the special case when $\Re(s)>1$, by applying (\ref{e0}), (\ref{e1}) and (\ref{e1'}) to (\ref{hurw}), we have
\begin{eqnarray*}
 \zeta(s,a)&=&\frac{1}{s-1}\left(\frac{\partial}{\partial z}\right)^{1-s}\!\!\left.\frac{ze^{az}}{e^z-1}\right|_{z=0}
 =\sum^\infty_{k=0}\left(\frac{\partial}{\partial z}\right)^{1-s}\!\!\left.\frac{ze^{(a+k)z}}{1-s}\right|_{z=0}\\
 &=&\sum^\infty_{k=0}\frac{\partial}{\partial\xi}e^{(a+k)\frac{\partial}{\partial\xi}}\left.\frac{\xi^{1-s}e^{\xi{z}}}{1-s}\right|_{\xi=0,z=0}\\
 &=&\sum^\infty_{k=0}\left(z+\frac{\partial}{\partial\xi}\right)e^{(a+k)\left(z+\frac{\partial}{\partial\xi}\right)}\left.\frac{\xi^{1-s}}{1-s}\right|_{\xi=0,z=0}\\
 &=&\sum^\infty_{k=0}\left(z+\frac{\partial}{\partial\xi}\right)\left.\frac{e^{(a+k)z}(\xi+a+k)^{1-s}}{1-s}\right|_{\xi=0,z=0}\\
 &=&\sum^\infty_{k=0}\left.\left(z\frac{e^{(a+k)z}(a+k)^{1-s}}{1-s}+e^{(a+k)z}(a+k)^{-s}\right)\right|_{z=0}\\
 &=&\sum^\infty_{k=0}\frac{1}{(a+k)^s}\quad(\Re(s)>1,\,a\not\in\mathbb{Z}_0^-).
\end{eqnarray*}
The result of this calculation shows that our definition is reasonable.

Similarly,
\begin{equation}\label{huzm}
 \zeta(s,a)=\sum^{m-1}_{n=0}\frac{1}{(n+a)^s}+\frac{1}{s-1}\left(\frac{\partial}{\partial{z}}\right)^{1-s}\!\!
 \left.\frac{ze^{(a+m)z}}{e^z-1}\right|_{z=0}\quad(m\in\mathbb{N}_0).
\end{equation}

\textbf{Definition 2.3.} The Riemann Zeta function $\zeta(s)$ can be defined in the complete form (This definition is valid for all complex $s$):
\begin{equation}\label{rize}
 \zeta(s):=\frac{1}{s-1}\left(\frac{\partial}{\partial z}\right)^{1-s}\!\!\left.\frac{ze^z}{e^z-1}\right|_{z=0}\quad(s\in\Omega:=\{s\in\mathbb{C}|s\neq1\}).
\end{equation}

\textbf{Definition 2.4.} Let $\zeta_z(s,a)$ be an analytic function defined by
\begin{equation}\label{rize'}
 \zeta_z(s,a):=\frac{1}{s-1}\left(\frac{\partial}{\partial z}\right)^{1-s}\!\!\frac{ze^{az}}{e^z-1}\quad(s\in\Omega:=\{s\in\mathbb{C}|s\neq1\},\,a\not\in\mathbb{Z}_0^-,\,z\in\mathbb{R}^1).
\end{equation}

In the special case when $\Re(s)>2$ and $z\leq0$, we have
\begin{equation}\label{rize"}
 \zeta_z(s,a)=\frac{1}{1-s}\sum^\infty_{n=0}\frac{ze^{(a+n)z}}{(a+n)^{s-1}}+\sum^\infty_{n=0}\frac{e^{(a+n)z}}{(a+n)^s}
 \quad(\Re(s)>2,\,a\not\in\mathbb{Z}_0^-,\,z\leq0).
\end{equation}

\textbf{Definition 2.5.} Let $\chi$ be a Dirichlet character modulo $q$. The Dirichlet L-function $L(s,\chi)$ can be defined in the complete form:
\begin{equation}\label{lfud}
 L(s,\chi):=\frac{1}{s-1}\left(\frac{\partial}{\partial{z}}\right)^{1-s}\sum^q_{k=1}\left.\chi(k)\frac{ze^{kz}}{e^{qz}-1}\right|_{z=0}
 \quad(s\in\mathbb{C}).
\end{equation}

This definition is valid for all complex $s$. In the special case when $\Re(s)>1$, by applying (\ref{e0}), (\ref{e1}) and (\ref{e1'}) to (\ref{lfud}), we have
\begin{eqnarray*}
 L(s,\chi)&=&\frac{1}{s-1}\left(\frac{\partial}{\partial{z}}\right)^{1-s}\sum^q_{k=1}\left.\chi(k)\frac{ze^{kz}}{e^{qz}-1}\right|_{z=0}\\
 &=&\sum^q_{k=1}\chi(k)\sum^\infty_{n=0}\left(\frac{\partial}{\partial z}\right)^{1-s}\!\!\left.\frac{ze^{(k+qn)z}}{1-s}\right|_{z=0}\\
 &=&\sum^q_{k=1}\chi(k)\sum^\infty_{n=0}\frac{\partial}{\partial\xi}e^{(k+qn)\frac{\partial}{\partial\xi}}\left.\frac{\xi^{1-s}e^{\xi{z}}}{1-s}\right|_{\xi=0,z=0}\\
 &=&\sum^q_{k=1}\chi(k)\sum^\infty_{n=0}\left(z+\frac{\partial}{\partial\xi}\right)
 e^{(k+qn)\left(z+\frac{\partial}{\partial\xi}\right)}\left.\frac{\xi^{1-s}}{1-s}\right|_{\xi=0,z=0}\\
 &=&\sum^q_{k=1}\chi(k)\sum^\infty_{n=0}\left(z+\frac{\partial}{\partial\xi}\right)\left.\frac{e^{(k+qn)z}(\xi+k+qn)^{1-s}}{1-s}\right|_{\xi=0,z=0}\\
 &=&\sum^q_{k=1}\chi(k)\sum^\infty_{n=0}\left.\left(z\frac{e^{(k+qn)z}(k+qn)^{1-s}}{1-s}+e^{(k+qn)z}(k+qn)^{-s}\right)\right|_{z=0}\\
 &=&\sum^\infty_{n=0}\sum^q_{k=1}\frac{\chi(k)}{(k+qn)^s}=\sum^\infty_{n=1}\sum^q_{k=1}\frac{\chi(k)}{(k+q(n-1))^s}\\
 &=&\sum^\infty_{n=1}\sum^{qn}_{k=1+q(n-1)}\chi(k-qn+q)\frac{1}{k^{s}}=\sum^\infty_{n=1}\frac{\chi(n)}{n^s}\quad(\Re(s)>1).
\end{eqnarray*}
The result of this calculation shows that our definition is reasonable.

Similarly,
\begin{equation}\label{lfjz}
 L(s,\chi)=\sum^{qm}_{n=1}\frac{\chi(n)}{n^s}+\frac{1}{s-1}\left(\frac{\partial}{\partial{z}}\right)^{1-s}
 \sum^q_{k=1}\left.\chi(k)\frac{ze^{(k+qm)z}}{e^{qz}-1}\right|_{z=0}\quad(m\in\mathbb{N}_0).
\end{equation}

\textbf{Definition 2.6.} Let $L_z(s,\chi)$ be an analytic function defined by
\begin{equation}\label{rild}
 L_z(s,\chi):=\frac{1}{s-1}\left(\frac{\partial}{\partial{z}}\right)^{1-s}\sum^q_{k=1}\chi(k)\frac{ze^{kz}}{e^{qz}-1}
 \quad(s\in\Omega:=\{s\in\mathbb{C}|s\neq1\},\,z\in\mathbb{R}^1).
\end{equation}

In the special case when $\Re(s)>2$ and $z\leq0$, since $\chi(k+q)=\chi(k)$, we have
\begin{eqnarray}\label{rild"}
 L_z(s,\chi) &=& \sum^q_{k=1}\chi(k)\left(\frac{1}{1-s}\sum^\infty_{n=0}\frac{ze^{(k+qn)z}}{(k+qn)^{s-1}}+\sum^\infty_{n=0}\frac{e^{(k+qn)z}}{(k+qn)^s}\right)\nonumber\\
 &=& \frac{z}{1-s}\sum^\infty_{n=1}\frac{\chi(n)}{n^{s-1}}e^{nz}+\sum^\infty_{n=1}\frac{\chi(n)}{n^s}e^{nz}.
\end{eqnarray}

\section{Main results}
\setcounter{equation}{0}
\renewcommand\theequation{3.\arabic{equation}}

\subsection{Alternating Zeta function}
\noindent

\textbf{Theorem 3.1.} For $r\in\mathbb{N}$
\begin{eqnarray}\label{d7.2a}
  \eta(2r+1) &=& \frac{2\times3^{2r}}{3^{2r+1}-1}\sum^r_{k=1}(-1)^{k-1}\frac{1}{(2k)!}\left(\frac{2\pi}{3}\right)^{2k}\eta(2r+1-2k)\nonumber\\
  & & +\,\frac{(2\pi)^{2r}}{3^{2r+1}-1}\sum^\infty_{k=1}(-1)^{r+k-1}\frac{E_{2k-1}(1)}{(2r+2k)!}\left(\frac{2\pi}{3}\right)^{2k}.
\end{eqnarray}
\begin{eqnarray}\label{d7.2b}
  \eta(2r+1) &=& \frac{2^{2r+1}}{2^{2r+1}-1}\sum^r_{k=1}(-1)^{k-1}\frac{1}{(2k)!}\left(\frac{\pi}{2}\right)^{2k}\eta(2r+1-2k)\nonumber\\
  & & +\,\frac{\pi^{2r}}{2^{2r+1}-1}\sum^\infty_{k=1}(-1)^{r+k-1}\frac{E_{2k-1}(1)}{(2r+2k)!}\left(\frac{\pi}{2}\right)^{2k}.
\end{eqnarray}
\begin{eqnarray}\label{d7.2c}
  \eta(2r+1) &=& \frac{3^{2r}(2^{2r}-1)}{3^{2r}(2^{2r-1}-1)+2^{2r-1}}\left[\sum^r_{k=1}\frac{(-1)^{k-1}}{(2k)!}\left(\frac{\pi}{3}\right)^{2k}\eta(2r+1-2k)\right.\nonumber\\
  & & \left.+\,\sum^\infty_{k=1}(-1)^{r+k-1}\frac{E_{2k-1}(1)}{2(2r+2k)!}\left(\frac{\pi}{3}\right)^{2r+2k}\right].
\end{eqnarray}

\textbf{Proof.} Letting $x=2c/3$ in (\ref{d7.1}), and by using (\ref{yc2.7}) we have (\ref{d7.2a}). Similarly we have (\ref{d7.2b}) and (\ref{d7.2c}) by using Lemma 2.5 and Lemma 2.7.

We can use $\eta(s)=(1-2^{1-s})\zeta(s)$ and (\ref{d7.2c}) to deduce the following theorem, which gives $\zeta(2r+1)$ recursively in terms of $\zeta(3),\zeta(5),\ldots,\zeta(2r-1)$:

\textbf{Theorem 3.2.} For $r\in\mathbb{N}$
\begin{eqnarray}\label{dy}
 \zeta(2r+1)&=& \frac{(-1)^r}{3^{2r}(2^{2r-1}-1)+2^{2r-1}}\left[\sum^{r-1}_{k=1}(-1)^{k-1}\frac{(2\pi)^{2r-2k}}{(2r-2k)!}3^{2k}(2^{2k}-1)\zeta(2k+1)\right.\nonumber\\
 & & \left.-\frac{(2\pi)^{2r}}{(2r)!}\ln2+(2\pi)^{2r}\sum^\infty_{k=1}(-1)^{k-1}\frac{E_{2k-1}(1)}{2(2r+2k)!}\left(\frac{\pi}{3}\right)^{2k}\right],
\end{eqnarray}
which provides a significantly simpler (and much more rapidly convergent) version of a main result of Dancs and He \cite[p. 194, Theorem 3.1]{Danh}:
\begin{eqnarray}\label{dh}
  & & (1-2^{-2r})\zeta(2r+1)\nonumber\\
  &=&\sum^{r-1}_{k=1}(-1)^k\frac{\pi^{2k}}{(2k+1)!}(2^{2k-2r}-1)\zeta(2r+1-2k)-\frac{(-1)^r\pi^{2r}\ln2}{(2r+1)!}\nonumber\\
  & & +\,(-1)^r\pi^{2r}\sum^\infty_{k=1}(-1)^{k-1}\frac{E_{2k-1}(1)}{2(2r+2k+1)!}\pi^{2k}\quad (r\in\mathbb{N}).
\end{eqnarray}

Similarly we have
\begin{eqnarray}\label{dyj}
 \zeta(2r+1)&=& \frac{(-1)^r2^{2r+1}}{(2^{2r}-1)(2^{2r+1}-1)}\left[\sum^{r-1}_{k=1}(-1)^{k-1}\frac{\pi^{2r-2k}}{(2r-2k)!}(2^{2k}-1)\zeta(2k+1)\right.\nonumber\\
 & & \left.-\frac{\pi^{2r}}{(2r)!}\ln2+\pi^{2r}\sum^\infty_{k=1}(-1)^{k-1}\frac{E_{2k-1}(1)}{2(2r+2k)!}\left(\frac{\pi}{2}\right)^{2k}\right]\quad (r\in\mathbb{N}),
\end{eqnarray}
which is a slightly corrected version of a main result proven in a different way by Tsumura \cite[p. 383, Theorem B]{Tsu}.

Since $E_{2k}(1)\equiv0$ for $k\in\mathbb{N}$,  we can use (\ref{d7.1}) and (\ref{44L}) to deduce the following theorem:

\textbf{Theorem 3.3.} For $m\in\mathbb{N}\setminus\{1\}$, in the interval $[-c,c]$ we have the following Fourier series expressions ([$\cdot$] is Gauss mark):
\begin{eqnarray}\label{44"}
 \sum^{\infty}_{n=1}(-1)^{n-1}\frac{1}{n^m}\cos\frac{n\pi{x}}{c}
 &=&\sum^{[m/2]}_{k=0}(-1)^{k}\frac{1}{(2k)!}\left(\frac{\pi{x}}{c}\right)^{2k}\eta(m-2k)\nonumber\\
 & &+\sum^\infty_{k=[m/2]+1}(-1)^k\frac{E_{2k-m}(1)}{2(2k)!}\left(\frac{\pi{x}}{c}\right)^{2k}.
\end{eqnarray}

Letting $x=c/3$ in (\ref{44"}), and by using (\ref{yc1.7}), we obtain

\textbf{Theorem 3.4.} For $m\in\mathbb{N}\setminus\{1\}$
\begin{eqnarray}\label{44"yc1.7}
  \eta(m) &=& \frac{3^{m-1}(2^{m-1}-1)}{3^{m-1}(2^{m-2}-1)+2^{m-2}}\left[\sum^{[m/2]}_{k=1}\frac{(-1)^{k-1}}{(2k)!}\left(\frac{\pi}{3}\right)^{2k}\eta(m-2k)\right.\nonumber\\
  & & \left.+\,\sum^\infty_{k=[m/2]+1}(-1)^{k-1}\frac{E_{2k-m}(1)}{2(2k)!}\left(\frac{\pi}{3}\right)^{2k}\right],
\end{eqnarray}
which gives $\eta(m)$ recursively in terms of $\eta(2),\eta(3),\eta(4),\ldots,\eta(m)$.
Since $E_{2k}(1)=0$ for $k\in\mathbb{N}$, letting $m=2,3,4,\cdots$ in (\ref{44"yc1.7}), we have
\begin{eqnarray*}
   \eta(2) &=& \frac{\pi^2}{12}, \\
   \eta(3) &=& \frac{3\pi^2}{22}\ln2+\frac{3\pi^2}{22}\sum^\infty_{k=1}(-1)^k\frac{E_{2k-1}(1)}{(2k+2)!}\left(\frac{\pi}{3}\right)^{2k},\\
   \eta(4) &=& \frac{189}{85}\left[\frac{\pi^2}{2!\cdot3^2}\eta(2)-\frac{\pi^4}{4!\cdot3^4}\frac{1}{2}\right]=\frac{7\pi^4}{720},\\
   &\cdots&,
\end{eqnarray*}
recursively. Since $E_{2k-1}(1)\neq0$, this provides some perspective on the difficulty of evaluating $\eta(2r+1)$ as opposed to $\eta(2r)$ for $r\in\mathbb{N}$.

By using (\ref{sch}), when $s=-2n,\,n\in\mathbb{N}$ in (\ref{yrdy}), we obtain
\[ \eta(-2n)=\left.\frac{\partial^{2n}}{\partial z^{2n}}\frac{e^z}{1+e^z}\right|_{z=0}=0.\]
When $s=-2n+1,\,n\in\mathbb{N}$ in (\ref{yrdy}), we obtain
\[ \eta(-2n+1)=\left.\frac{\partial^{2n-1}}{\partial z^{2n-1}}\frac{e^z}{1+e^z}\right|_{z=0}=\frac{1}{2}E_{2n-1}(1)=(2^{2n}-1)\frac{B_{2n}}{2n}.\]
Here $B_m$ are the Bernoulli numbers defined by the generating functions:
\begin{equation}\label{bnls}
 \frac{z}{e^z-1}=\sum^\infty_{m=0}B_m\frac{z^m}{m!}\quad(|z|<2\pi).
\end{equation}

Thus, we obtain the following asymptotic expansion for $\eta(s)$:

\textbf{Theorem 3.5.} For $s\in\mathbb{C}$
\begin{eqnarray}\label{yr1}
 \eta(s)\sim\sum^m_{n=1}\frac{(-1)^{n-1}}{n^s}
 +\frac{(-1)^m}{2}\left[\frac{1}{m^s}+\sum^\infty_{k=1}{-s\choose{2k-1}}\frac{E_{2k-1}(1)}{m^{s+2k-1}}\right]\quad(m\rightarrow\infty).
\end{eqnarray}

Let $-Li_s(-e^z)$ be an analytic function defined by
\begin{equation}\label{yrdy'}
-Li_s(-e^z):=\left(\frac{\partial}{\partial z}\right)^{-s}\!\!\frac{e^z}{1+e^z},\quad z\in\mathbb{R}^1,\,s\in\mathbb{C}.
\end{equation}
Then $\eta(s)=\left.-Li_s(-e^z)\right|_{z=0}$ and
\[\left.\frac{\partial^k}{\partial z^k}(-Li_s(-e^z))\right|_{z=0}=\eta(s-k)\quad(k\in\mathbb{N}).\]
Thus for $k\in\mathbb{N}$
\[\left.\frac{\partial^{2k-1}}{\partial x^{2k-1}}\cos\left(\frac{\pi x}{c}\frac{\partial}{\partial z}\right)(-Li_s(-e^z))\right|_{z=0,x=0}=0,\]
\[\left.\frac{\partial^{2k}}{\partial x^{2k}}\cos\left(\frac{\pi x}{c}\frac{\partial}{\partial z}\right)(-Li_s(-e^z))\right|_{z=0,x=0}=(-1)^k\left(\frac{\pi}{c}\right)^{2k}\eta(s-2k).\]
Here $c>0$ is a known real constant.

Therefore, we have the following Taylor expansion in the neighborhood of $x=0$:
\begin{equation}\label{yr2}
 \left.\cos\left(\frac{\pi x}{c}\frac{\partial}{\partial z}\right)(-Li_s(-e^z))\right|_{z=0}=\sum^\infty_{k=0}\frac{(-1)^k}{(2k)!}\left(\frac{\pi x}{c}\right)^{2k}\eta(s-2k)\quad (s\in\mathbb{C}).
\end{equation}
Let $x\in\mathbb{R}^1,\,c'>0$ be a real constant. By using Theorem 1.4, we have ($|x|\leq c'$)
\begin{equation}\label{yr3}
 \left.\cos\left(\frac{\pi x}{c}\frac{\partial}{\partial z}\right)(-Li_s(-e^z))\right|_{z=0}=\sum^\infty_{n=1}\frac{(-1)^{n-1}}{n^s}\cos\frac{n\pi x}{c}\quad (\Re(s)>1).
\end{equation}
Since $\frac{1}{2}E_n(1)=\eta(-n)$, by Theorem 3.3, we have ($x\in\mathbb{R}^1$ with $|x|\leq c$)
\begin{equation}\label{yr4}
 \sum^\infty_{n=1}\frac{(-1)^{n-1}}{n^m}\cos\frac{n\pi x}{c}=\sum^\infty_{k=0}\frac{(-1)^k}{(2k)!}\left(\frac{\pi x}{c}\right)^{2k}\eta(m-2k)\quad(m\in\mathbb{N}\setminus\{1\}).
\end{equation}
Letting $s=m$ in (\ref{yr3}), and by using formula (\ref{yr4}), we get $c'=c$. Thus we obtain the following generalization of Theorem 3.3:

\textbf{Theorem 3.6.} For $\Re(s)>1$ and $|x|\leq c\,(x\in\mathbb{R}^1)$, we have
\begin{equation}\label{yr5}
 \sum^\infty_{n=1}\frac{(-1)^{n-1}}{n^s}\cos\frac{n\pi x}{c}=\sum^\infty_{k=0}\frac{(-1)^k}{(2k)!}\left(\frac{\pi x}{c}\right)^{2k}\eta(s-2k).
\end{equation}

Similarly, let $\Phi(-1,s,a)$ be an analytic function defined by
\begin{equation}\label{yrdy"}
\Phi(-1,s,a):=\left(\frac{\partial}{\partial z}\right)^{-s}\!\!\left.\frac{e^{az}}{1+e^z}\right|_{z=0},\quad a\not\in\mathbb{Z}_0^-,\, s\in\mathbb{C}.
\end{equation}

\textbf{Theorem 3.7.} For $\Re(s)>1$ and $|x|\leq c\,(x\in\mathbb{R}^1)$, we have
\begin{equation}\label{yr6}
 \sum^\infty_{n=0}\frac{(-1)^n}{(n+a)^s}\cos\frac{(n+a)\pi x}{c}=\sum^\infty_{k=0}\frac{(-1)^k}{(2k)!}\left(\frac{\pi x}{c}\right)^{2k}\Phi(-1,s-2k,a).
\end{equation}
Here $a\neq0,-1,-2,\cdots$, and $\Phi(z,s,a)$ is the Lerch transcendent.

In the case when $a=1$, since $\Phi(-1,s,1)=\eta(s)$, Theorem 3.7 gives Theorem 3.6.

We can use Lemma 2.5 and Theorem 3.6 to deduce the following generalizations of the Theorem 3.1:

\textbf{Theorem 3.8.} For $s\in\mathbb{C}$

\begin{equation}\label{yrc1}
 \frac{3}{2}(1-3^{-s})\eta(s)=\sum^\infty_{k=1}\frac{(-1)^{k-1}}{(2k)!}\left(\frac{2\pi}{3}\right)^{2k}\eta(s-2k).
\end{equation}
\begin{equation}\label{yrc2}
 (1-2^{-s})\eta(s)=\sum^\infty_{k=1}\frac{(-1)^{k-1}}{(2k)!}\left(\frac{\pi}{2}\right)^{2k}\eta(s-2k).
\end{equation}
\begin{equation}\label{yrc3}
 \frac{1}{2}(1-2^{2-s}+3^{1-s})\zeta(s)=\sum^\infty_{k=1}\frac{(-1)^{k-1}}{(2k)!}\left(\frac{\pi}{3}\right)^{2k}\eta(s-2k).
\end{equation}

Letting $a=1/2$ in (\ref{yr6}), since $\Phi(-1,s,1/2)=2^s\beta(s)$, we obtain the following theorem:

\textbf{Theorem 3.9.} For $\Re(s)>0$ and $|x|\leq c/2\,(x\in\mathbb{R}^1)$, we have
\begin{equation}\label{yr7}
 \sum^\infty_{n=0}\frac{(-1)^n}{(2n+1)^s}\cos\frac{(2n+1)\pi x}{c}=\sum^\infty_{k=0}\frac{(-1)^k}{(2k)!}\left(\frac{\pi x}{c}\right)^{2k}\beta(s-2k).
\end{equation}
Here $\beta(s)$ is the Dirichlet Beta function, which is a holomorphic function on the whole complex plane. For $\Re(s)>0$, it is given by
\[\beta(s)=\sum^\infty_{n=0}\frac{(-1)^n}{(2n+1)^s}\quad(\Re(s)>0).\]

Letting $s=-n$ in (\ref{yrdy"}), and by using (\ref{sch}), we obtain
\begin{equation}\label{sch'}
 \Phi(-1,-n,a)=\frac{1}{2}E_n(a)\quad(n\in\mathbb{N}),
\end{equation}
where $E_n(x)$ are the Euler polynomials. In particular, $\eta(-n)=\frac{1}{2}E_n(1)$.

We can use (\ref{yrdy"}) and (\ref{sch'}) to deduce the following asymptotic expansion for $\Phi(-1,s,a)$:

\textbf{Theorem 3.10.} For $s\in\mathbb{C},\,a\not\in\mathbb{Z}_0^-$
\begin{eqnarray}\label{yr1a}
 \Phi(-1,s,a)\sim\sum^{m-1}_{n=0}\frac{(-1)^n}{(n+a)^s}
 +\frac{(-1)^m}{2}\left[\frac{1}{m^s}+\sum^\infty_{k=1}{-s\choose k}\frac{E_k(a)}{m^{s+k}}\right]\quad(m\rightarrow\infty).
\end{eqnarray}

When $a=1/2$ in (\ref{yr1a}), since $\Phi(-1,s,1/2)=2^s\beta(s)$, $E_k(1/2)=E_k/2^k$ and $E_{2k-1}=0$ for $k\in\mathbb{N}$, we obtain the following asymptotic expansion for $\beta(s)$:

\textbf{Theorem 3.11.} For $s\in\mathbb{C}$
\begin{eqnarray}\label{yr1b}
 \beta(s)\sim\sum^{m-1}_{n=0}\frac{(-1)^n}{(2n+1)^s}
 +\frac{(-1)^m}{2}\left[\frac{1}{(2m)^s}+\sum^\infty_{k=1}{-s\choose{2k}}\frac{E_{2k}}{(2m)^{s+2k}}\right]\quad(m\rightarrow\infty).
\end{eqnarray}
Here $E_n$ are the Euler numbers defined by the generating functions:
\begin{equation}\label{els}
 \frac{1}{\cosh z}=\frac{2}{e^z+e^{-z}}=\sum^\infty_{n=0}E_n\frac{z^n}{n!},\quad|z|<\frac{\pi}{2}.
\end{equation}

We can use (\ref{yr6}) and (\ref{sch'}) to deduce the following theorem:

\textbf{Theorem 3.12.} For $r\in\mathbb{N},\,a\not\in\mathbb{Z}_0^-$ and $|x|\leq c\,(x\in\mathbb{R}^1)$, we have
\begin{eqnarray}\label{yr8}
 \sum^\infty_{n=0}\frac{(-1)^n}{(n+a)^{2r+1}}\cos\frac{(n+a)\pi{x}}{c}
 &=& \sum^r_{k=0}\frac{(-1)^k}{(2k)!}\left(\frac{\pi{x}}{c}\right)^{2k}\Phi(-1,2r+1-2k,a)\nonumber\\
 & & +\sum^\infty_{k=1}(-1)^{r+k}\frac{E_{2k-1}(a)}{2(2r+2k)!}\left(\frac{\pi{x}}{c}\right)^{2r+2k};
\end{eqnarray}

\begin{eqnarray}\label{yr8'}
 \sum^\infty_{n=0}\frac{(-1)^n}{(n+a)^{2r}}\cos\frac{(n+a)\pi{x}}{c}
 &=& \sum^r_{k=0}\frac{(-1)^k}{(2k)!}\left(\frac{\pi{x}}{c}\right)^{2k}\Phi(-1,2r-2k,a)\nonumber\\
 & & +\sum^\infty_{k=1}(-1)^{r+k}\frac{E_{2k}(a)}{2(2r+2k)!}\left(\frac{\pi{x}}{c}\right)^{2r+2k}.
\end{eqnarray}

When $a=1$ in (\ref{yr8}) and (\ref{yr8'}), since $\Phi(-1,s,1)=\eta(s)$ and $E_{2k}(1)=0$, we obtain the Lemma 2.7 and Lemma 2.8 respectively.

Letting $a=1/2$ in (\ref{yr8}) and (\ref{yr8'}), since $\Phi(-1,s,1/2)=2^s\beta(s)$, $E_k(1/2)=E_k/2^k$ and $E_{2k-1}=0$ for $k\in\mathbb{N}$, we obtain the following results for $\beta(s)$:

\textbf{Theorem 3.13.} For $r\in\mathbb{N}$ and $|x|\leq c/2\,(x\in\mathbb{R}^1)$, we have
\begin{equation}\label{yr8b}
 \sum^\infty_{n=0}\frac{(-1)^n}{(2n+1)^{2r+1}}\cos\frac{(2n+1)\pi{x}}{c}
 =\sum^r_{k=0}\frac{(-1)^k}{(2k)!}\left(\frac{\pi{x}}{c}\right)^{2k}\beta(2r+1-2k);
\end{equation}

\begin{eqnarray}\label{yr8'b}
 \sum^\infty_{n=0}\frac{(-1)^n}{(2n+1)^{2r}}\cos\frac{(2n+1)\pi{x}}{c}
 &=& \sum^r_{k=0}\frac{(-1)^k}{(2k)!}\left(\frac{\pi{x}}{c}\right)^{2k}\beta(2r-2k)\nonumber\\
 & & +\sum^\infty_{k=1}\frac{(-1)^{r+k}E_{2k}}{2(2r+2k)!}\left(\frac{\pi{x}}{c}\right)^{2r+2k}.
\end{eqnarray}

When $x=c/2$, the series relations (\ref{yr8b}) and (\ref{yr8'b}) immediately yields

\begin{equation}\label{yr8d}
 \beta(2r+1)=\sum^r_{k=1}\frac{(-1)^{k-1}}{(2k)!}\left(\frac{\pi}{2}\right)^{2k}\beta(2r+1-2k)\quad(r\in\mathbb{N}),
\end{equation}
which gives $\beta(2r+1)$ recursively in terms of $\beta(3),\beta(5),\beta(7),\ldots,\beta(2r+1)$ using the initial value $\beta(1)=\pi/4$;
\begin{equation}\label{yr8'd}
 \beta(2r)=\sum^r_{k=1}\frac{(-1)^{k-1}}{(2k)!}\left(\frac{\pi}{2}\right)^{2k}\beta(2r-2k)
 +\sum^\infty_{k=1}\frac{(-1)^{r+k-1}E_{2k}}{2(2r+2k)!}\left(\frac{\pi}{2}\right)^{2r+2k},
\end{equation}
which gives $\beta(2r)$ recursively in terms of $\beta(2),\beta(4),\beta(6),\ldots,\beta(2r)$ using the initial value $\beta(0)=1/2$. Since $E_{2k}\neq0$ in (\ref{yr8'd}), this provides some perspective on the difficulty of evaluating $\beta(2r)$ as opposed to $\beta(2r+1)$.

Letting $s=-n$ in (\ref{yr1a}), since $\Phi(-1,-n,a)=E_n(a)/2$ and
\begin{equation}\label{elxz}
 E_n(x+y)=\sum^n_{k=0}{n\choose{k}}E_k(x)y^{n-k}\quad(n\in\mathbb{N}_0:=\mathbb{N}\cup\{0\}),
\end{equation}
we obtain the following summation formula:

\textbf{Theorem 3.14.} Let $a$ be any real number.  For $m\in\mathbb{N}$ and $n\in\mathbb{N}_0$
\begin{equation}\label{bfh}
 \sum^{m-1}_{k=0}(-1)^k(k+a)^n=\frac{1}{2}(E_n(a)-(-1)^mE_n(a+m)).
\end{equation}

\subsection{Hurwitz Zeta function}
\noindent

When $s=-n,\,n\in\mathbb{N}_0$ in (\ref{hurw}), we obtain
\begin{equation}\label{hurn}
 \zeta(-n,a)=-\frac{1}{n+1}\frac{\partial^{n+1}}{\partial z^{n+1}}\!\left.\frac{ze^{az}}{e^z-1}\right|_{z=0}=-\frac{1}{n+1}B_{n+1}(a),
\end{equation}
where $B_n(x)$ are the Bernoulli polynomials defined by the generating functions:
\begin{equation}\label{bnld}
  \frac{ze^{xz}}{e^z-1}=\sum^\infty_{n=0}B_n(x)\frac{z^n}{n!}\quad(|z|<2\pi).
\end{equation}
It is easily seen from the definition (\ref{hurw}) that
\begin{equation}\label{djd1}
 \lim_{s\rightarrow1}[(s-1)\zeta(s,a)]=\lim_{s\rightarrow1}\left(\frac{\partial}{\partial z}\right)^{1-s}\!\!\left.\frac{ze^{az}}{e^z-1}\right|_{z=0}
 =\left.\frac{ze^{az}}{e^z-1}\right|_{z=0}=B_0(a)=1.
\end{equation}
In other words, $\zeta(s,a)$ has a simple pole at $s=1$, and $B_0(a)=1$ is the residue of $\zeta(s,a)$ at the simple pole $s=1$.

Similarly, $\zeta(s)$ has a simple pole at $s=1$, and $B_0(1)=1$ is the residue of $\zeta(s)$ at the simple pole $s=1$.

We can use (\ref{huzm}) and (\ref{bnld}) to deduce the following asymptotic expansion for $\zeta(s,a)$:

\textbf{Theorem 3.15.} For $a\not\in\mathbb{Z}_0^-$, $s\in\mathbb{C}$ with $s\neq1$
\begin{eqnarray}\label{hurj}
 \zeta(s,a)\sim\sum^{m-1}_{n=0}\frac{1}{(n+a)^s}
 +\frac{1}{s-1}\sum^\infty_{k=0}{1-s\choose k}\frac{B_k(a)}{m^{s+k-1}}\quad(m\rightarrow\infty).
\end{eqnarray}

Since $B_0(1)=1,\,B_1(1)=1/2$, and for $k\in\mathbb{N}$, $B_{2k+1}(1)=0,\,B_{2k}(1)=B_{2k}$, when $a=1$, the Theorem 3.15 immediately yields
\begin{equation}\label{zezk}
 \zeta(s)\sim\sum^{m-1}_{n=1}\frac{1}{n^s}+\frac{m^{1-s}}{s-1}+\frac{m^{-s}}{2}
 +\frac{1}{s-1}\sum^\infty_{k=1}{1-s\choose2k}\frac{B_{2k}}{m^{s+2k-1}}\quad(m\rightarrow\infty),
\end{equation}
which is a slightly modified version of the following asymptotic expansion for $\zeta(s)$:
\begin{equation}\label{zezk'}
 \zeta(s)\sim\sum^{N-1}_{n=1}\frac{1}{n^s}+\frac{N^{1-s}}{s-1}+\frac{N^{-s}}{2}
 +N^{-s}\sum^\infty_{m=1}\frac{B_{2m}s^{\overline{2m-1}}}{(2m)!N^{2m-1}}\quad(N\rightarrow\infty),
\end{equation}
where $B_{2m}$ are the Bernoulli numbers and $s^{\overline{2m-1}}$ is a rising factorial. This expansion is valid for all complex $s$ and is often used to compute the Zeta function by using a large enough value of $N$.

\textbf{Theorem 3.16.} Let $a$ be any real number.  For $m\in\mathbb{N}$ and $n\in\mathbb{N}_0$
\begin{equation}\label{bfh"}
 \sum^{m-1}_{k=0}(k+a)^n=\frac{1}{n+1}(B_{n+1}(a+m)-B_{n+1}(a)).
\end{equation}

\textbf{Proof.} Letting $s=-n$ in (\ref{hurj}), it can be proved since $\zeta(-n,a)=-\frac{1}{n+1}B_{n+1}(a)$ and
\begin{equation}\label{bnhs}
 B_n(x+y)=\sum^n_{k=0}{n\choose{k}}B_k(x)y^{n-k}\quad(n\in\mathbb{N}_0).
\end{equation}

\textbf{Theorem 3.17.} For $a\not\in\mathbb{Z}_0^-$, $s\in\mathbb{C}$ with $s\neq1$
\begin{eqnarray}\label{hurj'}
 \zeta(s,a)\!\!&\sim&\!\!\sum^{m-1}_{n=0}\frac{1}{(n+a)^s}+\frac{(m+a)^{1-s}}{s-1}+\frac{(m+a)^{-s}}{2}\nonumber\\
            & & +\,\frac{1}{s-1}\sum^\infty_{k=1}{1-s\choose2k}\frac{B_{2k}}{(m+a)^{s+2k-1}}\quad(m\rightarrow\infty).
\end{eqnarray}

\textbf{Proof.} By applying (\ref{e1}) and (\ref{bnls}) to (\ref{huzm}), we have
\begin{eqnarray*}
 \zeta(s,a)&=&\sum^{m-1}_{n=0}\frac{1}{(n+a)^s}+\frac{1}{s-1}\left(m+a+\frac{\partial}{\partial z}\right)^{1-s}\!\!\left.\frac{z}{e^z-1}\right|_{z=0}\\
 &\!\!\sim\!\!&\sum^{m-1}_{n=0}\frac{1}{(n+a)^s}+\frac{(m+a)^{1-s}}{s-1}\sum^\infty_{k=0}{1-s\choose{k}}\frac{1}{(m+a)^k}
 \frac{\partial^k}{\partial{z}^k}\left.\frac{z}{e^z-1}\right|_{z=0}\\
 &=&\sum^{m-1}_{n=0}\frac{1}{(n+a)^s}+\frac{(m+a)^{1-s}}{s-1}\sum^\infty_{k=0}{1-s\choose{k}}\frac{B_k}{(m+a)^k}\quad(m\rightarrow\infty).
\end{eqnarray*}
Since $B_0=1,\,B_1=-1/2$, $B_{2k+1}=0$ for $k\in\mathbb{N}$, the Theorem 3.17 is proved.

Let $s^{(n)}=s^{\overline{n}}:=s(s+1)(s+2)\cdots(s+n-1)$ be the rising factorial, then
\[\frac{1}{s-1}{1-s\choose2k}=\frac{1}{s-1}{s-2+2k\choose2k}=\frac{1}{s-1}\frac{(s-1)^{\overline{2k}}}{(2k)!}=\frac{s^{\overline{2k-1}}}{(2k)!}\quad(s\neq1).\]
Thus the asymptotic expansion (\ref{hurj'}) can be written in the form:
\begin{eqnarray}\label{hurj"}
 \zeta(s,a)\!\!&\sim&\!\!\sum^{m-1}_{n=0}\frac{1}{(n+a)^s}+\frac{(m+a)^{1-s}}{s-1}+\frac{(m+a)^{-s}}{2}\nonumber\\
            & & +\,\sum^\infty_{k=1}\frac{B_{2k}\,s^{\overline{2k-1}}}{(2k)!\,(m+a)^{s+2k-1}}\quad(m\rightarrow\infty).
\end{eqnarray}

Letting $a=1$ in (\ref{hurj"}) with $m$ replaced by $N-1$, since $ \zeta(s,1)= \zeta(s)$, we obtain the asymptotic expansion (\ref{zezk'}) once again.

Similarly,

\textbf{Theorem 3.18.} For $s\in\mathbb{C},\,a\not\in\mathbb{Z}_0^-$
\begin{eqnarray}\label{yrsj}
 \Phi(-1,s,a)\!\!&\sim&\!\!\sum^{m-1}_{n=0}\frac{(-1)^n}{(n+a)^s}+\frac{(-1)^m}{2(m+a)^s}\nonumber\\
  & & +\,\frac{(-1)^m}{2}\sum^\infty_{k=1}\frac{E_{2k-1}(1)\,s^{\overline{2k-1}}}{(2k-1)!\,(m+a)^{s+2k-1}}\quad(m\rightarrow\infty).
\end{eqnarray}

By Definition 2.2 and Definition 2.4, we have $\left.\zeta_z(s,a)\right|_{z=0}=\zeta(s,a)$ and
\begin{equation}\label{wfze}
 \left.\frac{\partial^k}{\partial z^k}\zeta_z(s,a)\right|_{z=0}=\frac{s-1-k}{s-1}\zeta(s-k,a)\quad(k\in\mathbb{N},\,a\not\in\mathbb{Z}_0^-,\,s\neq1).
\end{equation}

From the differential relation (\ref{wfze}), we have for $k\in\mathbb{N},\,s\neq1$
\[ \left.\frac{\partial^{2k-1}}{\partial x^{2k-1}}\cos\left(\frac{\pi x}{c}\frac{\partial}{\partial z}\right)\zeta_z(s,a)\right|_{z=0,x=0}=0;\]
\[ \left.\frac{\partial^{2k}}{\partial x^{2k}}\cos\left(\frac{\pi{x}}{c}\frac{\partial}{\partial{z}}\right)\zeta_z(s,a)\right|_{z=0,x=0}
=(-1)^k\left(\frac{\pi}{c}\right)^{2k}\frac{s-1-2k}{s-1}\zeta(s-2k,a).\]
Here $c>0$ is a given real number. Therefore, we have the following Taylor expansion in the neighborhood of $x=0$:
\begin{equation}\label{taze}
 \left.\cos\left(\frac{\pi{x}}{c}\frac{\partial}{\partial{z}}\right)\zeta_z(s,a)\right|_{z=0}=\sum^\infty_{k=0}\frac{(-1)^k}{(2k)!}\left(\frac{\pi{x}}{c}\right)^{2k}
 \frac{s-1-2k}{s-1}\zeta(s-2k,a)\quad(s\neq1).
\end{equation}
On the other hand, by using (\ref{rize"}) and Theorem 1.4, since
\[\left.\cos\left(\frac{\pi{x}}{c}\frac{\partial}{\partial{z}}\right)e^{(a+n)z}\right|_{z=0}=\cos\frac{(a+n)\pi x}{c},\]
\[\left.\cos\left(\frac{\pi{x}}{c}\frac{\partial}{\partial{z}}\right)(ze^{(a+n)z})\right|_{z=0}
 =-\frac{\pi{x}}{c}\sin\frac{(a+n)\pi x}{c},\]
we obtain the following summation formula for $\Re(s)>2,\,x\in\Omega\subset\mathbb{R}^1$:
\begin{eqnarray}\label{sjqh}
 & &\left.\cos\left(\frac{\pi{x}}{c}\frac{\partial}{\partial{z}}\right)\zeta_z(s,a)\right|_{z=0}\nonumber\\
 &=&\sum^\infty_{n=0}\frac{1}{(a+n)^s}\cos\frac{(a+n)\pi x}{c}+\frac{\pi{x}/c}{s-1}\sum^\infty_{n=0}\frac{1}{(a+n)^{s-1}}\sin\frac{(a+n)\pi x}{c}.
\end{eqnarray}

By means of the  abstract operators $\cos(h\partial_x)$ and $\sin(h\partial_x)$, we can prove the following Katsurada's series representation:

Let $r\in\mathbb{N}$, $x\in\mathbb{R}^1$ with $-2c\leq{x}\leq2c$ ($c>0$ be a given real number). Then
\begin{eqnarray}\label{yl8}
 & & r\sum^\infty_{n=1}\frac{1}{n^{2r+1}}\cos\frac{n\pi{x}}{c}+\frac{\pi x}{2c}\sum^\infty_{n=1}\frac{1}{n^{2r}}\sin\frac{n\pi x}{c}\nonumber\\
 &=& \sum^{r-1}_{k=0}(-1)^k\,\frac{r-k}{(2k)!}\left(\frac{\pi{x}}{c}\right)^{2k}\zeta(2r+1-2k)\nonumber\\
 & & +\,(-1)^{r-1}\sum^\infty_{k=0}\frac{(2k)!}{(2r+2k)!}\frac{\zeta(2k)}{(2\pi)^{2k}}\left(\frac{\pi{x}}{c}\right)^{2r+2k},
\end{eqnarray}
which is a slightly modified version of a result proven in a significantly different way by Katsurada \cite[p. 81, Theorem 1]{Kats}.

Since
\[\zeta(2k)=(-1)^{k-1}\frac{(2\pi)^{2k}}{2(2k)!}B_{2k}\quad\mbox{and}\quad\zeta(-2k+1)=-\frac{1}{2k}B_{2k}\quad(k\in\mathbb{N}),\]
we can write the last term on the right-hand side of (\ref{yl8}) in the form:
\begin{eqnarray*}
  & & \sum^\infty_{k=r}(-1)^k\,\frac{r-k}{(2k)!}\left(\frac{\pi{x}}{c}\right)^{2k}\zeta(2r+1-2k) \\
  &=& \lim_{k\rightarrow0}a_k+\sum^\infty_{k=1}(-1)^{r+k}\,\frac{-k}{(2r+2k)!}\left(\frac{\pi{x}}{c}\right)^{2r+2k}\zeta(-2k+1) \\
  &=& \sum^\infty_{k=0}(-1)^{r+k}\,\frac{B_{2k}}{2(2r+2k)!}\left(\frac{\pi{x}}{c}\right)^{2r+2k} \\
  &=& (-1)^{r-1}\sum^\infty_{k=0}\frac{(2k)!}{(2r+2k)!}\frac{\zeta(2k)}{(2\pi)^{2k}}\left(\frac{\pi{x}}{c}\right)^{2r+2k},
\end{eqnarray*}
where
\[\lim_{k\rightarrow0}a_k=\lim_{k\rightarrow0}\left[(-1)^{r+k}\,\frac{-k}{(2r+2k)!}\left(\frac{\pi{x}}{c}\right)^{2r+2k}\zeta(-2k+1)\right]
=(-1)^r\,\frac{B_0}{2(2r)!}\left(\frac{\pi{x}}{c}\right)^{2r}.\]
Thus the Katsurada's formula (\ref{yl8}) can be written in the form:
\begin{eqnarray}\label{yl8g}
 & & r\sum^\infty_{n=1}\frac{1}{n^{2r+1}}\cos\frac{n\pi{x}}{c}+\frac{\pi x}{2c}\sum^\infty_{n=1}\frac{1}{n^{2r}}\sin\frac{n\pi x}{c}\nonumber\\
 &=& \sum^\infty_{k=0}(-1)^k\,\frac{r-k}{(2k)!}\left(\frac{\pi{x}}{c}\right)^{2k}\zeta(2r+1-2k)\quad(r\in\mathbb{N},\,|x|\leq2c).
\end{eqnarray}

By making use of (\ref{taze}), (\ref{sjqh}) and (\ref{yl8g}), we obtain the following theorem:

\textbf{Theorem 3.19.} Let $\Re(s)>1$. For $x\in\mathbb{R}^1$ with $|x|\leq2c$, we have
\begin{eqnarray}\label{zegh}
 & & s\sum^\infty_{n=0}\frac{1}{(a+n)^{s+1}}\cos\frac{(a+n)\pi x}{c}+\frac{\pi{x}}{c}\sum^\infty_{n=0}\frac{1}{(a+n)^s}\sin\frac{(a+n)\pi x}{c}\nonumber\\
 &=& \sum^\infty_{k=0}(-1)^k\frac{s-2k}{(2k)!}\left(\frac{\pi{x}}{c}\right)^{2k}\zeta(s+1-2k,a)\quad(\Re(s)>1).
\end{eqnarray}

By making use of (\ref{hurn}) and (\ref{zegh}), we obtain the following theorem:

\textbf{Theorem 3.20.} For $x\in\mathbb{R}^1$ with $|x|\leq2c$, we have
\begin{eqnarray}\label{zetl}
 & & r\sum^\infty_{n=0}\frac{1}{(a+n)^{2r+1}}\cos\frac{(a+n)\pi{x}}{c}+\frac{\pi x}{2c}\sum^\infty_{n=0}\frac{1}{(a+n)^{2r}}\sin\frac{(a+n)\pi x}{c}\nonumber\\
 &=& \sum^{r-1}_{k=0}(-1)^k\,\frac{r-k}{(2k)!}\left(\frac{\pi{x}}{c}\right)^{2k}\zeta(2r+1-2k,a)\nonumber\\
 & & +\,\sum^\infty_{k=0}(-1)^{r+k}\frac{B_{2k}(a)}{2(2r+2k)!}\left(\frac{\pi{x}}{c}\right)^{2r+2k}\quad(r\in\mathbb{N});
\end{eqnarray}
\begin{eqnarray}\label{zetl'}
 & & (2r-1)\sum^\infty_{n=0}\frac{1}{(a+n)^{2r}}\cos\frac{(a+n)\pi{x}}{c}+\frac{\pi{x}}{c}\sum^\infty_{n=0}\frac{1}{(a+n)^{2r-1}}\sin\frac{(a+n)\pi{x}}{c}\nonumber\\
 &=& \sum^{r-1}_{k=0}(-1)^k\,\frac{2r-1-2k}{(2k)!}\left(\frac{\pi{x}}{c}\right)^{2k}\zeta(2r-2k,a)\nonumber\\
 & & +\,\sum^\infty_{k=0}(-1)^{r+k}\frac{B_{2k+1}(a)}{(2r+2k)!}\left(\frac{\pi{x}}{c}\right)^{2r+2k}\quad(r\in\mathbb{N}\setminus\{1\}).
\end{eqnarray}

Since each series in (\ref{zegh}) is uniformly convergent with respect to $x$ on the closed interval $[-2c,2c]$ for $\Re(s)>2$, multiplying the both side of (\ref{zegh}) by $(\pi x/c)^{s-1}$, and executing termwise differentiation in them with respect to $x$, namely $\partial_x[(\pi x/c)^{s-1}\times(\ref{zegh})]$, then the sine series term is counteracted. Thus we can obtain the following theorem:

\textbf{Theorem 3.21.} Let $\Re(s)>2$. For $x\in\mathbb{R}^1$ with $|x|\leq2c$, we have
\begin{eqnarray}\label{zegh'}
 & &\sum^\infty_{n=0}\frac{s(s-1)}{(a+n)^{s+1}}\cos\frac{(a+n)\pi{x}}{c}
 +\left(\frac{\pi{x}}{c}\right)^2\sum^\infty_{n=0}\frac{1}{(a+n)^{s-1}}\cos\frac{(a+n)\pi{x}}{c}\nonumber\\
 &=&\sum^\infty_{k=0}(-1)^k\frac{(s-2k)(s+2k-1)}{(2k)!}\left(\frac{\pi{x}}{c}\right)^{2k}\zeta(s+1-2k,a)\quad(\Re(s)>2).
\end{eqnarray}

By making use of (\ref{hurn}) and (\ref{zegh'}), we obtain the following theorem:

\textbf{Theorem 3.22.} Let $r\in\mathbb{N}\setminus\{1\}$. For $x\in\mathbb{R}^1$ with $|x|\leq2c$, we have
\begin{eqnarray}\label{zetl"}
 & & \sum^\infty_{n=0}\frac{r(2r-1)}{(a+n)^{2r+1}}\cos\frac{(a+n)\pi{x}}{c}
 +\frac{1}{2}\left(\frac{\pi{x}}{c}\right)^2\sum^\infty_{n=0}\frac{1}{(a+n)^{2r-1}}\cos\frac{(a+n)\pi x}{c}\nonumber\\
 &=& \sum^{r-1}_{k=0}(-1)^k\,\frac{(r-k)(2r+2k-1)}{(2k)!}\left(\frac{\pi{x}}{c}\right)^{2k}\zeta(2r+1-2k,a)\nonumber\\
 & & +\,\sum^\infty_{k=0}(-1)^{r+k}\frac{(4r+2k-1)B_{2k}(a)}{2(2r+2k)!}\left(\frac{\pi{x}}{c}\right)^{2r+2k}.
\end{eqnarray}

\subsection{Dirichlet L-function}
\noindent

When $s=-n,\,n\in\mathbb{N}_0$ in (\ref{lfud}), we obtain
\begin{equation}\label{hurn"}
 L(-n,\chi)=-\frac{1}{n+1}\frac{\partial^{n+1}}{\partial z^{n+1}}\sum^q_{k=1}\left.\chi(k)\frac{ze^{kz}}{e^{qz}-1}\right|_{z=0}=-\frac{1}{n+1}B_{n+1,\chi},
\end{equation}
where $B_{n,\chi}$ are the generalized Bernoulli numbers defined by the generating functions:
\begin{equation}\label{bnld"}
 \sum^q_{k=1}\chi(k)\frac{ze^{kz}}{e^{qz}-1}=\sum^\infty_{n=0}B_{n,\chi}\frac{z^n}{n!}\quad(|z|<2\pi/q).
\end{equation}
It is easily seen from the definition (\ref{lfud}) that
\begin{equation}\label{djd2}
 \lim_{s\rightarrow1}[(s-1)L(s,\chi)]=\sum^q_{k=1}\left.\chi(k)\frac{ze^{kz}}{e^{qz}-1}\right|_{z=0}=B_{0,\chi}.
\end{equation}
Therefore, if $\chi$ is principal, then the corresponding Dirichlet L-function $L(s,\chi)$ has a simple pole at $s=1$, and $B_{0,\chi}\neq0$ is the residue of $L(s,\chi)$ at the simple pole $s=1$.

By applying (\ref{e1}) and (\ref{bnld"}) to (\ref{lfjz}), since
\[\frac{1}{s-1}{1-s\choose{k}}=\frac{(-1)^k}{s-1}{s-2+k\choose{k}}=\frac{(-1)^k(s-1)^{\overline{k}}}{k!\,(s-1)}=(-1)^k\frac{s^{\overline{k-1}}}{k!}\quad(s\neq1),\]
we obtain the following asymptotic expansion for $L(s,\chi)$:

\textbf{Theorem 3.23.} For $s\in\mathbb{C}$ with $s\neq1$
\begin{equation}\label{lfjz'}
 L(s,\chi)\sim\sum^{qm}_{n=1}\frac{\chi(n)}{n^s}+\frac{(qm)^{1-s}}{s-1}B_{0,\chi}
  +\sum^\infty_{k=1}(-1)^k\frac{B_{k,\chi}\,s^{\overline{k-1}}}{k!\,(qm)^{s+k-1}}\quad(m\rightarrow\infty).
\end{equation}

It's easy to see that (\ref{lfjz'}) can be written in the form:
\begin{equation}\label{lfjz"}
 L(s,\chi)\sim\sum^{qm}_{n=1}\frac{\chi(n)}{n^s}+\frac{(qm)^{1-s}}{s-1}
 \sum^\infty_{k=0}{1-s\choose{k}}\frac{B_{k,\chi}}{(qm)^k}\quad(m\rightarrow\infty).
\end{equation}

Letting $s=-n,\,n\in\mathbb{N}_0$ in (\ref{lfjz"}), since $L(-n,\chi)=-\frac{1}{n+1}B_{n+1,\chi}$, we obtain the following summation formula:

\textbf{Theorem 3.24.} For $n\in\mathbb{N}_0,\,m\in\mathbb{N}$
\begin{equation}\label{lfqh}
 \sum^{qm}_{k=1}\chi(k)k^n=\frac{1}{n+1}\sum^n_{k=0}{n+1\choose{k}}(qm)^{n+1-k}B_{k,\chi}.
\end{equation}

When $m=1$, we obtain the following recurrent formula of generalized Bernoulli numbers:
\begin{equation}\label{lfqh'}
 \frac{1}{n+1}\sum^n_{k=0}{n+1\choose{k}}q^{n+1-k}B_{k,\chi}=\sum^q_{k=1}\chi(k)k^n\quad(n\in\mathbb{N}_0),
\end{equation}
which gives $B_{n,\chi}$ recursively in terms of $B_{0,\chi},B_{1,\chi},B_{2,\chi},\ldots,B_{n,\chi}$.

By Definition 2.5 and Definition 2.6, we have $\left.L_z(s,\chi)\right|_{z=0}=L(s,\chi)$ and
\begin{equation}\label{wfld}
 \left.\frac{\partial^k}{\partial z^k}L_z(s,\chi)\right|_{z=0}=\frac{s-1-k}{s-1}L(s-k,\chi)\quad(k\in\mathbb{N},\,s\neq1).
\end{equation}

From the differential relation (\ref{wfld}), we have for $k\in\mathbb{N},\,s\neq1$
\[ \left.\frac{\partial^{2k-1}}{\partial x^{2k-1}}\cos\left(\frac{\pi x}{c}\frac{\partial}{\partial z}\right)L_z(s,\chi)\right|_{z=0,x=0}=0;\]
\[ \left.\frac{\partial^{2k}}{\partial x^{2k}}\cos\left(\frac{\pi{x}}{c}\frac{\partial}{\partial{z}}\right)L_z(s,\chi)\right|_{z=0,x=0}
=(-1)^k\left(\frac{\pi}{c}\right)^{2k}\frac{s-1-2k}{s-1}L(s-2k,\chi).\]
Here $c>0$ is a given real number. Therefore, we have the following Taylor expansion in the neighborhood of $x=0$:
\begin{equation}\label{tald}
 \left.\cos\left(\frac{\pi{x}}{c}\frac{\partial}{\partial{z}}\right)L_z(s,\chi)\right|_{z=0}=\sum^\infty_{k=0}\frac{(-1)^k}{(2k)!}\left(\frac{\pi{x}}{c}\right)^{2k}
 \frac{s-1-2k}{s-1}L(s-2k,\chi)\quad(s\neq1).
\end{equation}
Letting $s=2r+1,\,r\in\mathbb{N}$ in (\ref{tald}), since $L(-n, \chi)=-\frac{1}{n+1}B_{n+1,\chi},\,n\in\mathbb{N}_0$, we have
\begin{eqnarray}\label{ldcs}
 & & \left.\cos\left(\frac{\pi{x}}{c}\frac{\partial}{\partial{z}}\right)L_z(2r+1,\chi)\right|_{z=0}\nonumber\\
 &=& \frac{1}{r}\sum^{r-1}_{k=0}(-1)^k\,\frac{r-k}{(2k)!}\left(\frac{\pi{x}}{c}\right)^{2k}L(2r+1-2k,\chi)\nonumber\\
 & & +\,\frac{1}{2r}\sum^\infty_{k=0}(-1)^{r+k}\frac{B_{2k,\chi}}{(2r+2k)!}\left(\frac{\pi{x}}{c}\right)^{2r+2k}\quad(r\in\mathbb{N}).
\end{eqnarray}
By comparing (\ref{ldcs}) and (\ref{bnld"}), we have $|x|\leq2c/q$ in (\ref{ldcs}) for $r\in\mathbb{N}$, and thus $|x|\leq2c/q$ in (\ref{tald}).

On the other hand, by using (\ref{rild"}) and Theorem 1.4, we obtain the following summation formula for $\Re(s)>2,\,x\in\Omega\subset\mathbb{R}^1$:
\begin{equation}\label{sjld}
 \left.\cos\left(\frac{\pi{x}}{c}\frac{\partial}{\partial{z}}\right)L_z(s,\chi)\right|_{z=0}
 =\sum^\infty_{n=1}\frac{\chi(n)}{n^s}\cos\frac{n\pi x}{c}+\frac{\pi{x}/c}{s-1}\sum^\infty_{n=1}\frac{\chi(n)}{n^{s-1}}\sin\frac{n\pi x}{c}.
\end{equation}

By making use of (\ref{tald}) and (\ref{sjld}), we obtain the following theorem:

\textbf{Theorem 3.25.} Let $\Re(s)>1$. For $x\in\mathbb{R}^1$ with $|x|\leq2c/q$, we have
\begin{eqnarray}\label{ldgh}
 & & s\sum^\infty_{n=1}\frac{\chi(n)}{n^{s+1}}\cos\frac{n\pi x}{c}+\frac{\pi{x}}{c}\sum^\infty_{n=1}\frac{\chi(n)}{n^s}\sin\frac{n\pi x}{c}\nonumber\\
 &=& \sum^\infty_{k=0}(-1)^k\frac{s-2k}{(2k)!}\left(\frac{\pi{x}}{c}\right)^{2k}L(s+1-2k,\chi)\quad(\Re(s)>1).
\end{eqnarray}

By making use of (\ref{hurn"}) and (\ref{ldgh}), we obtain the following theorem:

\textbf{Theorem 3.26.} For $x\in\mathbb{R}^1$ with $|x|\leq2c/q$, we have
\begin{eqnarray}\label{ldtl}
 & & r\sum^\infty_{n=1}\frac{\chi(n)}{n^{2r+1}}\cos\frac{n\pi{x}}{c}+\frac{\pi x}{2c}\sum^\infty_{n=1}\frac{\chi(n)}{n^{2r}}\sin\frac{n\pi x}{c}\nonumber\\
 &=& \sum^{r-1}_{k=0}(-1)^k\,\frac{r-k}{(2k)!}\left(\frac{\pi{x}}{c}\right)^{2k}L(2r+1-2k,\chi)\nonumber\\
 & & +\,\sum^\infty_{k=0}(-1)^{r+k}\frac{B_{2k,\chi}}{2(2r+2k)!}\left(\frac{\pi{x}}{c}\right)^{2r+2k}\quad(r\in\mathbb{N});
\end{eqnarray}
\begin{eqnarray}\label{ldtl'}
 & & (2r-1)\sum^\infty_{n=1}\frac{\chi(n)}{n^{2r}}\cos\frac{n\pi{x}}{c}+\frac{\pi x}{c}\sum^\infty_{n=1}\frac{\chi(n)}{n^{2r-1}}\sin\frac{n\pi x}{c}\nonumber\\
 &=& \sum^{r-1}_{k=0}(-1)^k\,\frac{2r-1-2k}{(2k)!}\left(\frac{\pi{x}}{c}\right)^{2k}L(2r-2k,\chi)\nonumber\\
 & & +\,\sum^\infty_{k=0}(-1)^{r+k}\frac{B_{2k+1,\chi}}{(2r+2k)!}\left(\frac{\pi{x}}{c}\right)^{2r+2k}\quad(r\in\mathbb{N}\setminus\{1\}).
\end{eqnarray}

Since each series in (\ref{ldgh}) is uniformly convergent with respect to $x$ on the closed interval $[-2c/q,2c/q]$ for $\Re(s)>2$, multiplying the both side of (\ref{ldgh}) by $(\pi x/c)^{s-1}$, and executing termwise differentiation in them with respect to $x$, namely $\partial_x[(\pi x/c)^{s-1}\times(\ref{ldgh})]$, then the sine series term is counteracted. Thus we can obtain the following theorem:

\textbf{Theorem 3.27.} Let $\Re(s)>2$. For $x\in\mathbb{R}^1$ with $|x|\leq2c/q$, we have
\begin{eqnarray}\label{ldgh'}
 & &s(s-1)\sum^\infty_{n=1}\frac{\chi(n)}{n^{s+1}}\cos\frac{n\pi{x}}{c}
 +\left(\frac{\pi{x}}{c}\right)^2\sum^\infty_{n=1}\frac{\chi(n)}{n^{s-1}}\cos\frac{n\pi{x}}{c}\nonumber\\
 &=&\sum^\infty_{k=0}(-1)^k\frac{(s-2k)(s+2k-1)}{(2k)!}\left(\frac{\pi{x}}{c}\right)^{2k}L(s+1-2k,\chi).
\end{eqnarray}

By making use of (\ref{hurn"}) and (\ref{ldgh'}), we obtain the following theorem:

\textbf{Theorem 3.28.} For $x\in\mathbb{R}^1$ with $|x|\leq2c/q$, we have
\begin{eqnarray}\label{ldhu}
 & & r(2r-1)\sum^\infty_{n=1}\frac{\chi(n)}{n^{2r+1}}\cos\frac{n\pi{x}}{c}
 +\frac{1}{2}\left(\frac{\pi x}{c}\right)^2\sum^\infty_{n=1}\frac{\chi(n)}{n^{2r-1}}\cos\frac{n\pi x}{c}\nonumber\\
 &=& \sum^{r-1}_{k=0}(-1)^k\,\frac{(r-k)(2r+2k-1)}{(2k)!}\left(\frac{\pi{x}}{c}\right)^{2k}L(2r+1-2k,\chi)\nonumber\\
 & & +\,\sum^\infty_{k=0}(-1)^{r+k}\frac{(4r+2k-1)B_{2k,\chi}}{2(2r+2k)!}\left(\frac{\pi{x}}{c}\right)^{2r+2k}\quad(r\in\mathbb{N}\setminus\{1\});
\end{eqnarray}
\begin{eqnarray}\label{ldhu'}
 & & r(2r+1)\sum^\infty_{n=1}\frac{\chi(n)}{n^{2r+2}}\cos\frac{n\pi{x}}{c}
 +\frac{1}{2}\left(\frac{\pi x}{c}\right)^2\sum^\infty_{n=1}\frac{\chi(n)}{n^{2r}}\cos\frac{n\pi x}{c}\nonumber\\
 &=& \sum^r_{k=0}(-1)^k\,\frac{(2r+1-2k)(r+k)}{(2k)!}\left(\frac{\pi{x}}{c}\right)^{2k}L(2r+2-2k,\chi)\nonumber\\
 & & +\,\sum^\infty_{k=0}(-1)^{r+k+1}\frac{(2r+k+1)B_{2k+1,\chi}}{(2r+2k+2)!}\left(\frac{\pi{x}}{c}\right)^{2r+2k+2}\quad(r\in\mathbb{N}).
\end{eqnarray}

\subsection{Riemann Zeta function}
\noindent  

In the special case of Theorem 3.22 when $a=1$, since $\zeta(s,1)=\zeta(s)$, we obtain the following series representation:

\textbf{Theorem 3.29.} Let $r\in\mathbb{N}$, $x\in\mathbb{R}^1$ with $-2c\leq{x}<2c$ (or $-2c\leq x\leq2c$ for $r\in\mathbb{N}\setminus\{1\}$). Then
\begin{eqnarray}\label{yl15}
 & & r(2r-1)\sum^\infty_{n=1}\frac{1}{n^{2r+1}}\cos\frac{n\pi{x}}{c}
 +\frac{1}{2}\left(\frac{\pi x}{c}\right)^2\sum^\infty_{n=1}\frac{1}{n^{2r-1}}\cos\frac{n\pi x}{c}\nonumber\\
 &=& \sum^{r-1}_{k=0}(-1)^k\frac{(r-k)(2r+2k-1)}{(2k)!}\left(\frac{\pi{x}}{c}\right)^{2k}\zeta(2r+1-2k)\nonumber\\
 & & +\,(-1)^{r-1}\sum^\infty_{k=0}\frac{(2k)!(4r+2k-1)}{(2r+2k)!}\frac{\zeta(2k)}{(2\pi)^{2k}}\left(\frac{\pi{x}}{c}\right)^{2r+2k}.
\end{eqnarray}
Since $L(s,\chi)=\zeta(s)$ when $q=1$, we can also obtain the Theorem 3.29 from (\ref{ldhu}).

In particular, by making use of (\ref{yl15}) and Lemma 2.4, we can obtain the following important results:

\textbf{Theorem 3.30.} For $r\in\mathbb{N}\setminus\{1\}$
\begin{eqnarray}\label{yl21}
 \zeta(2r+1)&=& -\frac{2\pi^2(3^{2r-2}-1)}{r(2r-1)(3^{2r+1}-1)}\,\zeta(2r-1)+\frac{2\times3^{2r}}{r(2r-1)(3^{2r+1}-1)}\nonumber\\
 & & \times\left[\,\sum^{r-1}_{k=1}(-1)^{k-1}\frac{(r-k)(2r+2k-1)}{(2k)!}\left(\frac{2\pi}{3}\right)^{2k}\zeta(2r+1-2k)\right.\nonumber\\
 & & \left.+\,(-1)^r\left(\frac{2\pi}{3}\right)^{2r}\sum^\infty_{k=0}\frac{(2k)!(4r+2k-1)}{(2r+2k)!}\frac{\zeta(2k)}{3^{2k}}\right];
\end{eqnarray}
\begin{eqnarray}\label{yl22}
 \zeta(2r+1)&=& -\frac{(2^{2r-1}-2)\pi^2}{r(2r-1)(2^{4r+1}+2^{2r}-1)}\,\zeta(2r-1)+\frac{2^{4r+1}}{r(2r-1)(2^{4r+1}+2^{2r}-1)}\nonumber\\
 & & \times\left[\,\sum^{r-1}_{k=1}(-1)^{k-1}\frac{(r-k)(2r+2k-1)}{(2k)!}\left(\frac{\pi}{2}\right)^{2k}\zeta(2r+1-2k)\right.\nonumber\\
 & & \left.+\,(-1)^r\left(\frac{\pi}{2}\right)^{2r}\sum^\infty_{k=0}\frac{(2k)!(4r+2k-1)}{(2r+2k)!}\frac{\zeta(2k)}{4^{2k}}\right];
\end{eqnarray}
\begin{eqnarray}\label{yl23}
 \zeta(2r+1)&=& \frac{2\pi^2(6^{2r-2}-3^{2r-2}-2^{2r-2}+1)}{r(2r-1)(3^{2r}(2^{2r}+1)+2^{2r}-1)}\,\zeta(2r-1)\nonumber\\
 & & +\,\frac{2\times6^{2r}}{r(2r-1)(3^{2r}(2^{2r}+1)+2^{2r}-1)}\nonumber\\
 & & \times\left[\,\sum^{r-1}_{k=1}(-1)^{k-1}\frac{(r-k)(2r+2k-1)}{(2k)!}\left(\frac{\pi}{3}\right)^{2k}\zeta(2r+1-2k)\right.\nonumber\\
 & & \left.+\,(-1)^r\left(\frac{\pi}{3}\right)^{2r}\sum^\infty_{k=0}\frac{(2k)!(4r+2k-1)}{(2r+2k)!}\frac{\zeta(2k)}{6^{2k}}\right].
\end{eqnarray}

For example, letting $r=1,2,3,\cdots$ in (\ref{yl23}), we have
\begin{eqnarray*}
   \zeta(3) &=& -\frac{\pi^2}{6}\sum^\infty_{k=0}\frac{2k+3}{(2k+1)(2k+2)}\frac{\zeta(2k)}{6^{2k}}, \\
   \zeta(5) &=& \frac{8\pi^2}{87}\,\zeta(3)+\frac{\pi^4}{261}\sum^\infty_{k=0}\frac{(2k+7)\zeta(2k)}{(2k+1)(2k+2)\cdots(2k+4)6^{2k}},\\
   \zeta(7) &=& \frac{3124\pi^2}{29655}\,\zeta(5)-\frac{2\pi^4}{3295}\,\zeta(3)
    -\frac{16\pi^6}{88965}\sum^\infty_{k=0}\frac{(2k+11)\zeta(2k)}{(2k+1)(2k+2)\cdots(2k+6)6^{2k}},\\
   &\cdots&,
\end{eqnarray*}
recursively. Since $\zeta(2k)\rightarrow1$ as $k\rightarrow\infty$, each of these series representing $\zeta(2r+1)$ in (\ref{yl21}) to (\ref{yl23}) converges remarkably rapidly with its general term having the order estimate:
\[O(m^{-2k}\cdot k^{-2r+1})\qquad(k\rightarrow\infty;\quad m=3,4,6;\quad r\in\mathbb{N}).\]

Yan'an Second School, Yan'an 716000, Shaanxi, PR China

\emph{E-mail address}: guangqingbi@sohu.com

\end{CJK*}
\end{document}